\documentclass[12pt]{amsart}
\usepackage{graphicx}
\usepackage{verbatim}
\usepackage{textcomp}
\usepackage{amssymb}
\usepackage{cite}
\usepackage{color}
\usepackage[all]{xy}
\usepackage{graphicx}
\usepackage{color}
\usepackage{hyperref}
\hypersetup{
    colorlinks,
    citecolor=black,
    filecolor=black,
    linkcolor=black,
    urlcolor=black
}
\usepackage{tikz}
\usepackage[all]{xy}
\usepackage{graphicx}
\usetikzlibrary{backgrounds}
\newtheorem{theorem}{Theorem}[section]
\newtheorem{lemma}[theorem]{Lemma}
\newtheorem{proposition}[theorem]{Proposition}
\newtheorem{corollary}[theorem]{Corollary}   

\newtheorem{claim}[theorem]{Claim}
\allowdisplaybreaks
\newtheorem{notation}[theorem]{Notation}

\theoremstyle{definition}
\newtheorem{definition}[theorem]{Definition}

\theoremstyle{remark}
\newtheorem{remark}[theorem]{Remark}

\numberwithin{equation}{section}

\newtheorem*{theorem*}{Theorem}
{{\sc Proof of Lemma~\ref{tri1}.}}%
{{\qed} \\}
{{\sc Proof of Theorem~\ref{regularity}.}}%
{{\qed} \\}
{{\sc Proof of Theorem~\ref{main}.}}%
{{\qed} \\}
\newenvironment{proofof(i)}%
    {{\sc Proof of $(i)$.}}%
  {{\qed} \\}  
  
  \newenvironment{proofof(iv)}%
    {{\sc Proof of $(iv)$.}}%
  {{\qed} \\}  

\newcommand{\R}{\mathbb R}

\newcommand{\Sp}{\mathbb S}

\newcommand{\B}{\mathcal B}

\title{Regularity of harmonic maps from Polyhedra to CAT(1) Spaces}

\thanks{This work began as part of the workshop ``Women in Geometry" (15w5135) at the Banff International Research Station in November of 2016. We are grateful to BIRS for the opportunity to attend and for the excellent working environment. 
CB, CM were supported in part by NSF grants DMS-1308420 and DMS-1406332 respectively, and LH was supported by NSF grants DMS 1308837 and DMS 1452477. AF was supported in part by an NSERC Discovery Grant. PS was supported in part by an NSERC PGS D scholarship and a UBC Four Year Doctoral Fellowship. YZ was supported in part by an AWM-NSF Travel Grant. This material is also based upon work supported by NSF DMS-1440140 while CB and AF were in residence at the Mathematical Sciences Research Institute in Berkeley, California, during the Spring 2016 semester.}
\author[Breiner]{Christine Breiner}
\address{Department of Mathematics \\
                 Fordham University \\
                 Bronx, NY  10458}
\email{cbreiner@fordham.edu}
\author[Fraser]{Ailana Fraser}
\address{Department of Mathematics \\
                 University of British Columbia \\
                 Vancouver, BC V6T 1Z2}
\email{afraser@math.ubc.ca}
\author[Huang]{Lan-Hsuan Huang}
\address{Department of Mathematics\\ University of Connecticut\\ Storrs, CT 06269}
\email{lan-hsuan.huang@uconn.edu}
\author[Mese]{Chikako Mese}
\address{Johns Hopkins University\\
Department of Mathematics\\
3400 N. Charles Street\\
Baltimore, MD  21218}
\email{cmese@math.jhu.edu}
\author[Sargent]{Pam Sargent}
\address{Department of Mathematics \\
                 University of British Columbia \\
                 Vancouver, BC V6T 1Z2}
\email{psargent@math.ubc.ca}
\author[Zhang]{Yingying Zhang}
\address{Johns Hopkins University\\
Department of Mathematics\\
3400 N. Charles Street\\
Baltimore, MD  21218}
\email{yzhang@math.jhu.edu}

\begin{document}
\maketitle
\begin{abstract}
We determine regularity results for energy minimizing maps from an $n$-dimensional Riemannian polyhedral complex $X$ into a CAT(1) space. Provided that the metric on $X$ is Lipschitz regular, we prove H\"older regularity with  H\"older constant and exponent dependent on the total energy of the map and the metric on the domain. Moreover, at points away from the $(n-2)$-skeleton, we improve the regularity to locally Lipschitz. Finally, for points $x \in X^{(k)}$ with $k \leq n-2$, we demonstrate that the H\"older exponent depends on geometric and combinatorial data of the link of $x \in X$.
\end{abstract}
\section{Introduction}
A natural notion of energy for a map
between geometric  spaces  is defined by  measuring the total stretch of the map at each point of the domain and then integrating it over the domain.    Harmonic maps are critical points of the energy functional.  They
can be seen as both a generalization of harmonic functions in complex analysis and a higher dimensional analogue of parameterized geodesics in Riemannian geometry.   In the absence of a   totally geodesic map, a harmonic map is perhaps the most natural way to map one given geometric space  into another. 

The celebrated work of Eells and Sampson \cite{eells-sampson} initiated a wide interest in the study of  harmonic maps  between Riemannian manifolds, and harmonic maps have proven to be a useful tool in geometry.  A more recent development is the harmonic map theory for non-smooth spaces.  
The seminal works of   
Gromov-Schoen \cite{gromov-schoen} and  Korevaar-Schoen \cite{korevaar-schoen1}  consider harmonic maps from a Riemannian domain into a non-Riemannian  target.  Further exploration of harmonic map theory in the singular setting  includes works of Jost \cite{jost}, J. Chen \cite{chen}, Eells-Fuglede \cite{eells-fuglede} and Daskalopoulos-Mese \cite{daskal-meseCAG}. The above mentioned works all assume non-positivity of curvature (NPC) in the target space.  In this paper, the goal is to investigate the regularity issues of harmonic maps in the case when the target curvature is bounded above by a constant that is not necessarily 0.  In this direction, we mention earlier works of Serbinowski \cite{serbinowski} for harmonic maps from  Riemannian  manifold domains and Fuglede \cite{fuglede-TAMS,fuglede-CVPDE} for polyhedral domains.   

By understanding the regularity of harmonic maps, we can realize the  potential applications of harmonic map theory.   The key issue is to prove regularity theorems strong enough to be able to apply differential geometric methods. Applications of harmonic maps  already in the literature  include  those in  rigidity problems  (for example, \cite{siu}, \cite{corlette}, \cite{gromov-schoen}) and in Teichm\"{u}ller theory (for example, \cite{wolf}, \cite{donaldson}, \cite{daskal-meseHR}) amongst others.   Our goal is to apply harmonic map theory in a more general  setting (namely for CAT(1) targets) than the  NPC targets  considered in the above mentioned applications.   Indeed, in the follow-up  of this paper \cite{Paper2}, we prove a generalization to the metric space setting of Sacks and Uhlenbeck's celebrated work \cite{saks-uhlenbeck} on the bubbling phenomena for harmonic maps.  The generalization of Sacks and Uhlenbeck's work has important connections to the non-smooth uniformization problem (cf. \cite{bonk-kleiner} and references therein) which in turn is related to the Canon conjecture and the asymptotic geometry of negatively curved spaces.  Details of these connections are provided in the introduction of \cite{Paper2}.
We now state our main theorems.

\begin{theorem}[H\"older Regularity]\label{Holder}
 Let $B(r)$ be a ball of radius $r$ around a point $x$ in an admissible complex $X$ endowed with a Lipschitz Riemannian metric $g$ and let $(Y, d)$ be  a CAT(1) space and $\varrho \in (0,1)$. If $f:(B(r),g) \to \B_{\tau}(P) \subset Y$ is an energy minimizing map, where $0<\tau< \frac \pi 4$, then there
exist $C>0, \gamma>0$ such that
\[
 d(f(x),f(y)) \leq C|x-y|^\gamma \text{ for all } x, y \in B(\varrho r).
\]The constants $C, \gamma$ depend only on the total energy $E^f_g$ of the map, $(B(r),g)$ and $\varrho$.
\end{theorem}
\begin{remark}
Note that in the statement of the main theorems, the radius $r$ of the ball $B(r)$ is measured with respect to the Euclidean metric $\delta_{ij}$ on each cell. 
\end{remark}
Fuglede proves a similar result in \cite{fuglede-TAMS,fuglede-CVPDE}.   The main improvements of our results are the following:   First, the metric on the domain space is more general; more specifically, the metrics considered in this paper are only assumed to be Lipschitz continuous while  Fuglede considers  simplex-wise smooth metrics  (cf.~page 380, subsection ``Maps into metric spaces" in   \cite{fuglede-CVPDE}).  We hope that this will lead to wider applications for the theory of harmonic maps from polyhedral domains.  Second, and more importantly, we explicitly give the dependence of the H\"{o}lder constant and exponent on the total energy of the map.  This statement in the special case of NPC targets has been crucial in the applications of harmonic map theory.  In particular, the explicit dependence leads to a compactness result for a family of harmonic maps with uniformly bounded energy (see \cite[Lemma 3.1]{Paper2}).  Moreover, we can deduce the existence of tangent maps associated to harmonic maps (see Proposition \ref{tangentmapprop}).

We further remark that our proof uses  very different techniques from those in \cite{fuglede-TAMS,fuglede-CVPDE}.  Specifically,  we take advantage of the work done by Daskalopoulos and Mese for NPC targets 
in \cite{daskal-meseCAG},
using the order function and a Campanato type theorem to prove the H\"older regularity.   One of the advantages of this method is that, on high dimensional faces, we can improve the regularity to gain Lipschitz control, as given in Theorem \ref{LipThm}. Moreover, as in \cite{daskal-meseCAG}, for points in the lower dimensional skeleta, we provide a lower bound on the H\"older exponent of the minimizing map in terms of the first eigenvalue of the link of the normal strata of the skeleton, $\lambda_1^N$. 
\begin{theorem}[Lipschitz Regularity]\label{LipThm}
 Let $B(r)$ be a ball of radius $r$ around a point $x$ in an admissible complex $X$ endowed with a Lipschitz Riemannian metric $g$ and let $(Y,d)$ be a CAT(1) space. Suppose that $f:(B(r),g) \to \B_{\tau}(P) \subset Y$ is an energy minimizing map where $0< \tau < \frac \pi 4$. 
 \begin{enumerate}
 \item For $x \in X -X^{(n-2)}$, let $\bar d$ denote the distance of $x$ to $X^{(n-2)}$. Then for $\varrho \in (0,1)$ and $d' \leq \min\{\varrho r,  \varrho \bar d\}$, $f$ is Lipschitz continuous
 in $B(d')$ with Lipschitz constant depending on the total energy $E^f_g$ of the map $f$, $(B(r),g)$, and $d'$.
 \item For $x \in X^{(k)} - X^{(k-1)}$ and $k = 0, 1, \dots, n-2$, let $\bar d$ denote the distance of $x$ to $X^{(k-1)}$. Then for $\varrho \in (0,1)$ and $d' \leq \min\{ \varrho r,\varrho \bar d\}$, $f$ is H\"older continuous in $B(d')$ with H\"older exponent and constant depending on the total energy $E^f_g$ of the map $f$, $(B(r),g)$, and $d'$. More precisely, the H\"older exponent $\alpha$ has lower bound given by the following: If $\lambda_1^N \geq \beta \, (>\beta)$ then $\alpha\, (\alpha +n-k-2) \geq \beta \,(> \beta)$. In particular, if $\lambda_1^N \geq n-k-1$, then $f$ is Lipschitz continuous in a neighborhood of $x$.
 \end{enumerate}
 \end{theorem} 
To understand the second item, for $x \in X^{(k)} - X^{(k-1)}$, let $N=N(x)$ denote the link of $X^{(k)}$ at $x$ with metric induced by the Lipschitz Riemannian metric on $X$. Set
\[
\lambda_1^N:= \inf_{Q \in Y} \lambda_1(N,T_QY),
\]where $\lambda_1(N,T_QY)$ denotes the first eigenvalue of the Laplacian of $N$ with values in the tangent cone of $Y$ at $Q$. For more details, see section \ref{Lipss}.

Serbinowski \cite{serbinowski}, in an unpublished thesis, proves Lipschitz regularity from a Riemannian domain. Again, our proof is quite different from his.   Since the regularity theorems above are local results, we also obtain the following.

\begin{corollary}
We have the same conclusions of all of the previous theorems if we replace the assumption that  $Y$ is CAT(1) by $Y$ is locally CAT(1).
\end{corollary}

The paper is organized as follows. In section \ref{domain_target}, we define the domains and targets of interest and prove a few key estimates on CAT(1) spaces. Section \ref{energy_sec} includes background and necessary references for defining the energy and minimizing maps into metric spaces. This section also includes the definition of the cone over $Y$ and important distance relations. In section \ref{monotoness}, we prove a monotonicity formula for minimizing maps into CAT(1) spaces. In section \ref{Holderss}, we use the monotonicity formula to prove Theorem \ref{Holder}. Section \ref{Improved} uses Theorem \ref{Holder} to improve the monotonicity result which in turn allows us to improve the H\"older regularity so that the H\"older exponent is given by the order of the map. In section \ref{Tangentss}, we determine a tangent map construction using the cone over $Y$, where the existence of a tangent map is given by the H\"older regularity. Finally, in section \ref{Lipss}, we use the tangent map construction and the improved H\"older regularity to prove Theorem \ref{LipThm}.

\section{Domain and target spaces}
\label{domain_target}

\subsection{Admissible cell complexes and local models}

Throughout the paper, $X$ will denote an admissible $n$-dimensional  cell complex (i.e. a dimensionally homogeneously, locally $(n-1)$-chainable  convex cell complex) with a Lipschitz continuous Riemannian metric defined on each cell.  We refer to  \cite[Section 2.2]{daskal-meseCAG} for more details.     In particular, since the regularity theorems we prove are local, we will study harmonic maps from a ``local model"   that represents a neighborhood of a point of $X$.  We refer  the reader to \cite[Section 2.1]{daskal-meseCAG} for the precise formulation of a local model, but will briefly describe this here.  To do so,  we inductively define a  $k$-dimensional  \emph{conical cell}.   First, a $1$-dimensional conical cell is either the interval $[0,\infty)$ or the interval $(-\infty,0]$.  Having defined $(k-1)$-dimensional conical cells, we define a $k$-dimensional conical cell $C$ as a subset of ${\mathbb R}^k$ with the following properties:
\begin{itemize}
\item[(i)] The set $C$ is non-empty and closed.  
\item[(ii)]  The set $C$ is conical; i.e. if $x \in C$, then $tx \in C$ for $t \geq 0$.  
\item[(iii)]  The  intersection of $C$ with the unit sphere ${\mathbb S}^{k-1} \subset {\mathbb R}^k$  is geodesically convex (with respect the standard metric on ${\mathbb S}^{k-1}$).
\item[(iv)] The boundary $\partial C$  of $C$  is a finite union of $\{c_i\}$ where each  $c_i$ is  a subset of a $(k-1)$-dimensional hyperplane $H_i$ of ${\mathbb R}^k$ containing the origin such that if we identify $H_i \subset {\mathbb R}^k$ with ${\mathbb R}^{k-1}=\{(x^1, \dots, x^{k-1},0)\} \subset {\mathbb R}^k$, (via an orthogonal transformation which takes $H_i$ to ${\mathbb R}^{k-1}$), then  $c_i$ is a $(k-1)$-dimensional conical cell.  We will  say that $c_i$ is a  \emph{$(k-1)$-dimensional boundary cell of $C$}.
\end{itemize}
An  \emph{$l$-dimensional boundary cell}   of a $k$-dimensional conical cell $C$ is $H \cap C$, where $H$ is again a hyperplane of $\mathbb R^k$ containing the origin, such that there exists an orthogonal transformation of ${\mathbb R}^k$ which takes $H \cap C$ into ${\mathbb R}^l=\{(x^1, \dots, x^l,0,\dots,0)\} \subset {\mathbb R}^k$ but there exists no orthogonal transformation which takes $H \cap C$ into ${\mathbb R}^{l-1}=\{(x^1, \dots, x^{l-1},0,\dots,0)\} \subset {\mathbb R}^k$.  The union of $l$-dimensional boundary cells is called the \emph{$l$-skeleton} of $C$. 

Note that since $\partial C$ bounds a conical cell, the hyperplanes $H_i$ containing $c_i$ are linearly independent in the sense that the set of normal vectors defining the hyperplanes are all linearly independent. Indeed, one may consider $C$ as the intersection of appropriately oriented half-spaces, each with boundary one of the $H_i$. A $k$-dimensional conical cell $C$  is said to have \emph{codimension $\nu$} if $\partial C = \cup_{i=1}^\nu c_i$. In that case,
there exists  a hyperplane $H$ of ${\mathbb R}^k$ containing the origin and an orthogonal transformation $T$ of ${\mathbb R}^k$  such that $T(H \cap \partial C)$ is equal to ${\mathbb R}^{k-\nu} = \{(x^1, \dots, x^{k-\nu},0,\dots, 0) \}\subset {\mathbb R}^k$. We let $D:=T^{-1}(\mathbb R^{k-\nu})$.

A \emph{dimension-$n$, codimension-$\nu$ local model}  (of a neighborhood of a point in an $n$-dimensional cell complex) is ${\bf B}:=\bigsqcup W/\sim$, i.e.~ a disjoint union of a finite number ${\mathcal W} =\{W\}$ of $n$-dimensional conical cells of codimension $\nu$  modulo an equivalence relation $\sim$.   We   refer to $W \in {\mathcal W}$ as a \emph{wedge}.  The equivalence relation $\sim$ is defined  by a finite set of isometries  $\{\varphi\}$ where each $\varphi$ maps a    boundary cell  of one wedge to a boundary   cell of another wedge. Note that the equivalence relation implies that we may consider a single $D$ as belonging to the local model ${\bf B}$.

We assume ${\bf B}$ is \emph{admissible},  i.e.~ whenever $W \in \mathcal W$ and $S$ is a $(n-2)$-skeleton of $W$,  $\bigsqcup W \backslash S\slash \sim$ is connected.

Each wedge $W$ of ${\bf B}$ is a subset of ${\mathbb R}^n$ and therefore ${\bf B}$ comes equipped with the Euclidean metric (because each $W$ inherits the Euclidean metric  from ${\mathbb R}^n$). Let ${\bf B}(r)$ denote the ball of radius $r$, with respect to the Euclidean metric, centered at the origin of ${\bf B}$. 
Throughout the rest of the paper $B_x(\sigma)$ will denote a Euclidean ball in ${\bf B}$, centered at $x$ and of radius $\sigma$.  
Furthermore, using the coordinates inherited from ${\mathbb R}^n$, we can define a Riemannian metric $g$ on ${\bf B}$ by defining component functions $(g_{ij})$ on each wedge $W$.  We say $g$  is a \emph{Lipschitz Riemannian metric} on ${\bf B}$ if on each $W$
\[
|g_{ij}(x)-g_{ij}(\bar{x})|\leq c|x-\bar{x}|, \ \ \forall x, \bar{x} \in W.
\]
As explained in \cite[Proposition 2.1]{daskal-meseCAG}, we can and will often assume that the Lipschitz metrics are \emph{normalized}, i.e.
   \begin{equation} \label{lipbd}
  |g_{ij}(x)-\delta_{ij}| \leq c\sigma \text{ for } |x| \leq \sigma.
\end{equation}
Thus, for a \emph{normalized} Lipschitz metric,  \[
g_{ij}(0)=\delta_{ij}.
\]

Lastly, we say $\lambda \in (0,1]$ is an \emph{ellipticity constant} for $g$ if for each wedge $W$ and for $x \in W$, 
\[
\lambda |\xi|^2 \leq \sum_{i,j=1}^n g_{ij} (x)\xi^i \xi^j \leq \lambda^{-1} |\xi|^2.
\]

\subsection{CAT(1) spaces}

Given a complete metric space $(Y,d)$, $Y$ is called a \emph{geodesic space} if for each $P, Q \in Y$, there exists a curve $\gamma_{PQ}$ such that the length of $\gamma_{PQ}$ is exactly $d(P, Q)$.  We call $\gamma_{PQ}$ a \emph{geodesic} between $P$ and $Q$. 
\begin{remark}
For ease of notation, we will often denote $d(P,Q)$ by $d_{PQ}$. 
\end{remark}

We determine a weak notion of an upper sectional curvature bound on $Y$ by using comparison triangles. Given any three points $P,Q,R \in Y$ such that $d_{PQ}+d_{QR}+d_{RS} < 2\pi$, the \emph{geodesic triangle} $\triangle PQR$ is the triangle in $Y$ with sides given by the geodesics $\gamma_{PQ}, \gamma_{QR}, \gamma_{RS}$.  

Let $\triangle \tilde{P}\tilde{Q}\tilde{R}$ denote a geodesic triangle on the standard sphere $\Sp^2$ such that $d_{PQ}=d_{\tilde{P}\tilde{Q}}$, $d_{QR}=d_{\tilde{Q}\tilde{R}}$  and $d_{RP}=d_{\tilde{R}\tilde{P}}$. We call $\triangle \tilde{P}\tilde{Q}\tilde{R}$ a \emph{comparison triangle} for the geodesic triangle $\triangle PQR$. Note that a comparison triangle is convex since the perimeter of the geodesic triangle is less than $2\pi$.

\begin{definition}
Given a geodesic space $(Y,d)$ and a geodesic $\gamma_{PQ}$ with $d_{PQ}<\pi$, for $\tau \in [0,1]$ let $(1-\tau) P + \tau Q$ denote the point on $\gamma_{PQ}$ at distance $\tau d_{PQ}$ from $P$. That is,
\[
d((1-\tau) P + \tau Q,P) = \tau d_{PQ}.
\]
\end{definition}
\begin{definition}
Let $(Y,d)$ be a complete geodesic space. Then $Y$ is a CAT(1) space if:

Given any geodesic triangle $\triangle PQR$ (with perimeter less than $2\pi$) and a comparison triangle $\triangle \tilde{P}\tilde{Q}\tilde{R}$ in $\Sp^2$, 
\begin{equation} \label{cat}
d_{P_tR_s} \leq d_{\tilde{P}_t\tilde{R}_s}
\end{equation}
where 
\begin{eqnarray*}
P_t=(1-t)P+tQ,  & &  R_s=(1-s)R+sQ,\\
\tilde{P}_t=(1-t)\tilde{P}+t\tilde{Q}, & &  \tilde{R}_s=(1-s)\tilde{R}+s\tilde{Q}.
\end{eqnarray*}
A complete geodesic space $Y$ is said to be \emph{locally}   CAT(1) if every point of $Y$ has a geodesically convex CAT(1) neighborhood.
\end{definition}
We conclude this section with a few key estimates that we will use later in the paper. 
The first estimate appeared in the thesis of \cite[Estimate II]{serbinowski} without proof. See \cite[Lemma A.4]{Paper2} for a proof.

\begin{lemma}
 \label{lemma:estimateII}
Let $\triangle PQS$ be a geodesic triangle in a CAT(1) space $(Y,d)$. For a  pair of numbers $0\le \eta, \eta' \le 1$ define 
\begin{align*}
	P_{\eta'} &= (1-\eta') P + \eta' Q\\
	S_\eta& = (1-\eta) S + \eta Q. 
\end{align*}
Then 
\begin{align}\begin{split}\label{estII}
	d^2 (P_{\eta'}, S_{\eta}) &\le  \frac{\sin^2((1-\eta) d_{QS})}{\sin^2 d_{QS}} (d_{PS}^2- (d_{QS}-d_{QP})^2) +  \left((1-\eta)(d_{QS}-d_{QP}) + (\eta' - \eta) d_{QS}\right)^2  \\
	&\quad+ \mathrm{Cub}\left( d_{PS}, d_{QS}-d_{QP}, \eta-\eta' \right).
\end{split}\end{align}
\end{lemma}



As an immediate consequence, we have the following lemma.
\begin{lemma}\label{tri1}
Let $\triangle PQS$ be a geodesic triangle in a CAT(1) space $(Y,d)$. For $0\le \eta, \eta' \le 1$ and $P_{\eta'}, S_\eta$ as above,
\begin{align*}
	d^2 (P_{\eta'}, S_{\eta}) &\le(1 -2\eta + \eta d^2_{QS})d_{PS}^2  -2(\eta-\eta')(d_{QS}-d_{QP})d_{QS} + (\eta' - \eta)^2 d_{QS}^2  \\
	&\quad+\mathrm{Quad}(\eta, \eta') \mathrm{Quad}(d_{PS}, d_{QS}-d_{QP})+ \mathrm{Cub}\left( d_{PS}, d_{QS}-d_{QP}, \eta-\eta' \right).
\end{align*}
\end{lemma}
\begin{proof}
By Taylor expansion, $\sin((1-\eta)d_{QS}) = \sin d_{QS} - \eta d_{QS} \cos d_{QS} + O(\eta^2)$. Since $\frac{a}{\sin a} \geq 1$ and $\cos a \geq 1-\frac {a^2}2$ for $0 \leq a < \pi$, 
\[
 \frac{\sin^2((1-\eta) d_{QS})}{\sin^2 d_{QS}} = \left(1 - \eta \frac {d_{QS}}{\sin d_{QS}} \cos d_{QS} +O(\eta^2)\right)^2\leq 1- 2\eta + \eta d^2_{QS} + O(\eta^2). 
\]Substituting into \eqref{estII} implies that
 \begin{align*}
 	d^2 (P_{\eta'}, S_{\eta}) &\le  \left(1- 2\eta + \eta d^2_{QS}\right) (d_{PS}^2- (d_{QS}-d_{QP})^2) +  \left((1-\eta)(d_{QS}-d_{QP}) + (\eta' - \eta) d_{QS}\right)^2  \\
	&\quad + \eta^2\mathrm{Quad}(d_{PS}, d_{QS}-d_{QP})+ \mathrm{Cub}\left( d_{PS}, d_{QS}-d_{QP}, \eta-\eta' \right).
\end{align*}Expanding the quadratic term and collecting the remaining like terms implies the result.
\end{proof}

We conclude this section with a convexity bound.
\begin{lemma}\label{lemmacvx}
Let $\triangle PQR$ be a geodesic triangle in a CAT(1) space $(Y,d)$.  If $d_{PQ},  d_{PR}< \frac \pi 2$, then
\begin{equation} \label{fcvx}
\frac 1{8} \cos \left(d_{Q_{\frac 12}P}\right)d^2_{QR} \leq \frac 12 (d^2_{RP} + d^2_{QP})-d^2_{Q_{\frac 12}P}
\end{equation}
where $Q_{\frac 12}$ denotes the midpoint between $Q$ and $R$.
\end{lemma}
\begin{proof}
By the triangle comparison, it suffices to prove  inequality (\ref{fcvx})  assuming that $\triangle PQR$ is a geodesic triangle on the unit sphere. 
Let $\gamma(s)$ be an arclength parameterized geodesic on the sphere. Let
\[
d_{\gamma}(s):=d(\gamma(s),P)
\]
and assume that for all $s$, $d_\gamma(s)< \frac \pi 2$.
The function $d_{\gamma}(s)$ satisfies 
\[
(\cos d_{\gamma}(s)) '' = -\cos d_{\gamma}(s).
\]
Direct computation shows that
\begin{eqnarray*}
(d_{\gamma}(s))''& = &  \frac{1}{\tan d_{\gamma}(s)} (1-(d'_{\gamma}(s))^2)\\
(d^2_{\gamma}(s))'' & = & 2\frac{d_{\gamma}(s)}{\tan d_{\gamma}(s)} \left( 1-(d'_{\gamma}(s))^2\right) + 2(d'_{\gamma}(s))^2.
\end{eqnarray*}
Thus,
\[
(d^2_{\gamma}(s)+2\cos d_{\gamma}(s))'' = 2 \left( \frac{d_{\gamma}(s)}{\tan d_{\gamma}(s)} -\cos d_{\gamma}(s) \right) +2 \left( 1- \frac{d_{\gamma}(s)}{\tan d_{\gamma}(s)}  \right) (d'_{\gamma}(s))^2.
\]
Now let $\sigma(t)$ be a constant speed parameterization of the geodesic with $\sigma(0)=Q$ and $\sigma(1)=R$.   Thus, $\gamma(s):=\sigma(s/\delta)$, where $\delta=d_{QR}$, is an arclength parameterized geodesic.  With 
\[
d(t):=d(\sigma(t),P),
\]
  the chain rule implies that
\[
(d^2(t)+2\cos d(t)) ''=  2\delta^2 \left( \frac{d(t)}{\tan d(t)} -\cos d(t) \right)+ 2 \delta^2\left(1- \frac{d(t)}{\tan d(t)}  \right) ( d'(t))^2.
\]
Since $0\leq d(t) < \frac \pi 2$,
\[
\cos d(t) \leq   \frac{d(t)}{\tan d(t)} \leq 1.
 \]
Therefore
\[
(d^2(t)+2\cos d(t))''  \geq 0.
\]
The convexity of $t \mapsto d^2(t)+ 2\cos d(t)$ implies that
\begin{equation} \label{cvxfunny}
d^2_{Q_{\frac12 P}} + 2 \cos d_{Q_{\frac 12}P} \leq \frac{1}{2} \left( d^2_{QP}+2\cos d_{QP} + d^2_{RP}+ 2\cos d_{RP} \right).
\end{equation}

Using the identity on the sphere and a double angle formula,
\[
\cos d_{Q_{\frac 12}P} = \frac{\sin \frac{d_{QR}}{2}}{\sin d_{QR}} 
\cos d_{QP} +\frac{\sin \frac{d_{QR}}{2}}{\sin d_{QR}} 
\cos d_{RP}= \frac 1{2 \cos  \frac{d_{QR}}{2}}\left( \cos d_{QP} + \cos d_{RP}\right).
\]Since $\cos a \leq 1- \frac {a^2}{4}$ for $0 \leq a < \frac \pi 4$,

\[
\cos d_{QP} + \cos d_{RP} =2 \cos d_{Q_{\frac 12}P} \cos \frac{d_{QR}}2 \leq 2 \cos \left(d_{Q_{\frac 12}P}\right)- \frac 18 \cos \left(d_{Q_{\frac 12}P}\right)d^2_{QR}. 
\]The desired inequality follows from inserting the above into \eqref{cvxfunny}.
\end{proof}

\section{Sobolev space and the energy density}\label{energy_sec}
 In the seminal work of Korevaar-Schoen (cf. \cite[Chapter 1]{korevaar-schoen1}) the authors define the energy density and directional energies for maps from Riemannian manifolds into metric spaces. Using \cite[Proposition 2.1]{daskal-meseCAG}, these definitions immediately extend to include maps from an admissible complex $X$ (cf. \cite[Section 2]{daskal-meseCAG}). Following the usual convention, we say $f \in W^{1,2}(\Omega,Y)$ if $f \in L^2(\Omega)$ and the energy density is finite. We then write $|\nabla f|^2_g(x)$ in place
of the energy density function and let
\[
E^f_g =
\int_{{\bf B}(r)}|\nabla f|^2_g d\mu_g.
\]
 For a set $S \subset {\bf B}(r)$, let
\[
E^f_g[S]=\int_S |\nabla f|_g^2 d\mu_g.
\]

To study energy minimizing maps, we use the notion of the trace of $f$, for $f\in W^{1,2}(\Omega,Y)$, as defined in \cite{korevaar-schoen1} and \cite{eells-fuglede}. 
We denote the space of admissible maps $W^{1,2}_f(\Omega, \B):= \{ h \in W^{1,2}(\Omega, \B) : d(f,h) \in W^{1,2}_0(\Omega)\}$. 

\begin{definition}\label{mindef}Let $\Omega$ be a compact domain in an admissible complex with Lipschitz Riemannian metric $g$ and $(Y,d)$ be a CAT(1) space. 
A finite energy map $f:\Omega \rightarrow \B_\tau(P) \subset Y$ is \emph{energy minimizing} if $f$ minimizes energy amongst maps in
$W^{1,2}_f(\Omega,\overline{\B_{\tau}(P)})$. 
\end{definition}

The existence and uniqueness of energy minimizers from Riemannian domains appeared in the thesis \cite{serbinowski} and the same result from
Riemannian complexes into small balls in a CAT(1) space was established in \cite{fuglede-AASF}. We verify the existence and uniqueness in the Riemannian case in the appendix of \cite{Paper2}.

\begin{remark}
Note that unlike the definition in \cite{serbinowski}, the comparison maps in Definition \ref{mindef} not only have the same trace as $f$ but also map into the same ball. The reason that we define energy minimizing maps in this way is that, unlike in the NPC setting, the projection map onto convex domains in a CAT(1) space is not globally distance decreasing. Therefore, one cannot guarantee that a minimizer in the class $W^{1,2}_f(\Omega,Y)$ maps into the closure of $\B_\tau(P)$ without some extra hypotheses. For simplicity, we define a minimizer by considering only competitors in the smaller class of maps.
\end{remark}

\subsection{The pullback metric} 

The directional energies are defined in a fashion similar to the energy density function. See \cite{korevaar-schoen1} for the definition of the directional energy and of the pull-back inner product $\pi$ when $Y$ is an NPC space.

We use the triangle comparison in CAT(1) spaces to demonstrate that directional energies and the pull-back inner product are well defined for finite energy maps into CAT(1) spaces. The next lemma appeared in
 \cite[Lemma 3.6]{meseCAG} and is a consequence of \eqref{cat}. We include the proof here both for completeness and because we have simplified the proof.

\begin{lemma} \label{gb}
Let $Y$ be a CAT(1) space.  For every $\epsilon_0>0$, there exists  $\delta_0>0$  such that if  $P,Q,R,S \in Y$ with   $\max\{d_{PQ},d_{QR},d_{RS},d_{PS}\}\leq\delta_0$, then
\[
d^2_{PR} +d^2_{QS} \leq d^2_{PQ} +d^2_{QR} +d^2_{RS} +d^2_{PS} +\epsilon_0 \delta_0^2.
\]
\end{lemma}

\begin{proof}
 By comparing a geodesic quadrilateral $\square PQRS$ in $Y$ to a comparison quadrilateral $\square \tilde P \tilde Q \tilde R \tilde S$ in $\Sp^2$ (and noting \cite{reshetnyak} which says that the  pairwise distance of points on $\square PQRS$ is bounded by the distance of the corresponding pair in $\square \tilde P \tilde Q \tilde R \tilde S$), it is sufficient to prove the assertion when $Y=\Sp^2$.   Suppose the assertion is not true on $\Sp^2$.  Then there exists $\epsilon_0>0$ and a sequence $P_i, Q_i, R_i, S_i$ with $ \max\{d_{P_iQ_i},d_{Q_iR_i},d_{R_iS_i},d_{P_iS_i}\} \leq\delta_i\rightarrow 0$ such that
 \[
d^2_{P_iR_i} +d^2_{Q_iS_i} > d^2_{P_iQ_i} +d^2_{Q_iR_i} +d^2_{R_iS_i} +d^2_{P_iS_i} +\epsilon_0 \delta_i^2.
\]
For each $i$, denote by $\frac{1}{\delta_i}\Sp^2$ the rescaling of  the unit sphere $\Sp^2$ by a factor of  $\frac{1}{\delta_i}$ and let $P_i',Q_i',R_i',S_i' \in \frac{1}{\delta_i}\Sp^2$ be the corresponding points to $P_i,Q_i,R_i,S_i \in \Sp^2$ respectively.  Thus,
\[
d^2_{P_i'R_i'} +d^2_{Q_i'S_i'} >d^2_{P_i'Q_i'} +d^2_{Q_i'R_i'} +d^2_{R_i'S_i'} +d^2_{P_i'S_i'} +\epsilon_0.
\]
and $\max\{d_{P_i'Q_i'},d_{Q_i'R_i'},d_{R_i'S_i'},d_{P_i'S_i'}\}=1$.  This is a contradiction since the Gauss curvature of the sphere $\frac{1}{\delta_i}\Sp^2$   goes to 0  as $\delta_i\rightarrow 0$ and 
\[
d^2_{\bar{P}\bar{R}} +d^2_{\bar{Q}\bar{S}} \leq d^2_{\bar{P}\bar{Q}} +d^2_{\bar{Q}\bar{R}} +d^2_{\bar{R}\bar{S}} +d^2_{\bar{S}\bar{P}}
\]
for every $\bar{P}, \bar{Q}, \bar{R},\bar{S} \in \R^2$.
\end{proof}

\begin{lemma} \label{parallelogram}
Let  $f:({\bf B}(r),g) \to Y$ be a finite energy map and $(Y,d)$ a CAT(1) space. Then the parallelogram identity
\[
|f_*(Z + W)|^2_g + |f_*(Z - W)|^2_g = 2|f_*(Z)|^2_g + 2|f_*(W)|^2_g
\]
 holds for a.e. $x \in {\bf B}(r)$ and any pair of Lipschitz vector fields  $Z,W$ on ${\bf B}(r)$.
 \end{lemma}
 
 \begin{proof}
Fix $\epsilon_0>0$.
Let  $\epsilon \mapsto x_1(x,\epsilon)$, $\epsilon \mapsto  x_2(x,\epsilon)$ and $\epsilon \mapsto x_3(x,\epsilon)$ be the flow induced by the vector fields $Z$, $Z+W$ and $W$ 
respectively with $x_1(x,0)=x_2(x,0)=x_3(x,0)=x$. 
By \cite[Lemma 1.9.2]{korevaar-schoen1},  $\epsilon \mapsto f(x_i(x,\epsilon))$ is continuous at $\epsilon=0$,  $i=1,2,3$ for a.e. $x \in \Omega$.  
For such $x$,  apply Lemma~\ref{gb} with $P = f(x)$, $Q = f(x_1(x,\epsilon))$, $R = f(x_2(x,\epsilon))$, $S = f(x_3(x,\epsilon))$, divide by $\epsilon^2$ and multiply 
the resulting inequality by $\phi \in C^{\infty}_c(\Omega)$.  Now following the argument of \cite[Lemma 2.3.1]{korevaar-schoen1}, we conclude that
 \[
|f_*(Z + W)|^2_g(x) + |f_*(Z - W)|^2_g(x)  \leq  2|f_*(Z)|^2_g(x) + 2|f_*(W)|^2_g(x)+\epsilon_0 \Delta
\]
where 
\[
\Delta\geq\max \{|f_*(Z + W)|^2_g(x), |f_*(Z - W)|^2_g(x), |f_*(Z)|^2_g(x), |f_*(W)|^2_g(x)\}.
\]
  Since $\epsilon_0>0$ is arbitrary, we conclude that
\[
|f_*(Z + W)|^2_g + |f_*(Z - W)|^2_g  \leq  2|f_*(Z)|^2_g + 2|f_*(W)|^2_g.
\]
Repeat using $Z+W$ and $Z-W$ in place of $Z$ and $W$ to get the opposite
inequality. 
 \end{proof}

 \begin{lemma} \label{pullbackmetric}
 Let $f:({\bf B}(r),g) \to Y$ be a finite energy map and let $Z,W$ be Lipschitz vector fields on ${\bf B}(r)$. The operator $^g\pi^f$ defined by
 \[
 ^g\pi^f(Z,W) :=\frac{1}{2} |f_*(Z +W)|^2_g - \frac{1}{2}|f_*(Z -W)|^2_g
 \]
 is  symmetric, bilinear, non-negative and tensorial.
   \end{lemma}
 \begin{proof}
Using Lemma~\ref{parallelogram}, we can follow the  proof of \cite[Theorem 2.3.2]{korevaar-schoen1}.
 \end{proof}

\begin{notation} Let $\left\{ \frac{\partial}{\partial x^1}, ...,\frac{\partial}{\partial x^n}  \right\}$  be the standard Euclidean basis defined on each wedge inherited from ${\R}^n$ and $\delta$  the standard Euclidean metric.  Set
\[
\frac{\partial f}{\partial x^i} \cdot \frac{\partial f}{\partial
x^j}=\ ^\delta \pi^f\left(\frac{\partial}{\partial x^i},
\frac{\partial}{\partial x^j} \right) \ \mbox{ and } \ \left|
\frac{\partial f}{\partial x^i} \right|^2 = \frac{\partial
f}{\partial x^i} \cdot \frac{\partial f}{\partial x^i}.
\]
Similarly for the standard Euclidean polar coordinates $(r, \theta_1, \dots, \theta_{n-1})$ on each wedge we denote

\[
 \frac{\partial
f}{\partial x^k} \cdot \frac{\partial f}{\partial r}=\ ^{\delta} \pi^f \left( \frac{\partial
}{\partial x^k} ,\frac{\partial }{\partial r} \right),
 \
\left| \frac{\partial f}{\partial r} \right|^2= \frac{\partial
f}{\partial r} \cdot \frac{\partial f}{\partial r}=\ ^{\delta} \pi^f \left( \frac{\partial
}{\partial r} ,\frac{\partial }{\partial r} \right)
\]
and
\[
 \frac{\partial
f}{\partial {\theta_i}} \cdot \frac{\partial f}{\partial {\theta_j}}=\ ^\delta \pi^f\left(\frac{\partial}{\partial {\theta_i}},
\frac{\partial}{\partial {\theta_j}} \right).
\]
\end{notation}
Note that the energy density with respect to the metric $g$ is given by
\[
|\nabla f|^2_g
=\sum_{i,j} g^{ij} \frac{\partial f}{\partial x^i} \cdot
\frac{\partial f}{\partial x^j},
\]
whereas the energy density with respect to the Euclidean metric is given by
\[
|\nabla f|^2=|\nabla f|^2_{\delta}
=\sum_{i} \left| \frac{\partial f}{\partial x^i} \right|^2.
\]

\subsection{The cone over Y and energy comparisons}
\label{conesection} 

We denote by ${\mathcal C}Y$  the metric cone over $Y$. Topologically,
${\mathcal C}Y$ is defined by
\[
{\mathcal C}Y = Y \times [0,\infty)/Y \times \{0\}.
\]
A point in ${\mathcal C}Y$ is a pair $[P, t]$ for  $P \in Y$ and $t \in  [0,\infty)$, with $[P, 0]$
and $[Q, 0]$ representing the same point in ${\mathcal C}Y$ for all $P,Q \in Y$. We endow ${\mathcal C}Y$ with
a distance function $D$ defined by
$D^2([P, t], [Q, s]) = t^2 + s^2 - 2ts \cos \min(d_{PQ}, \pi)$.
It is well known that when $Y$ is a CAT(1) space, the metric space $({\mathcal C}Y,D)$ is an NPC space.

For $P,Q \in Y$ with $d_{PQ}<\pi/2$,  
\begin{equation}\label{D_d_compare}
\frac 12 \leq \frac{D^2([P,1],[Q,1])}{d^2_{PQ}} \leq 1,
\end{equation}
\begin{equation}\label{projectiondistance}
\lim_{P \to Q}\frac{D^2([P,1],[Q,1])}{d^2_{PQ}} = \lim_{P \to Q}\frac{2(1- \cos (d_{PQ}))}{d^2_{PQ}}=1,
\end{equation}
i.e., $Y$ is isometrically embedded into $Y \times \{1\}$ in an infinitesimal sense. 
Moreover, for $d_{PQ}$ small,
\begin{equation}\label{D_d_compare2}
d^2_{PQ}(1-d^2_{PQ}) \leq D^2([P,1],[Q,1]).
\end{equation}
\begin{definition}\label{lifteddef}
For any map $w:\Omega \to Y$, we let $\overline w:\Omega \to Y \times \{1\}$ be given by $\overline w(x) = [w(x),1]$. We call $\overline w$ the \emph{lifted map} of $w$.
\end{definition}If $w \in W^{1,2}(\Omega,Y)$ then $\overline w \in W^{1,2}(\Omega,\mathcal CY)$
and the definition of energy implies that
\begin{equation}\label{liftenergy}
{^d}E^w_g[\Omega] = {^D}E_g^{\overline w}[\Omega].
\end{equation}

We let $\Pi:\mathcal CY \to Y\times \{1\}$ denote the projection map $\Pi([P,t]) = [P,1]$. Then for $d_{PQ}< \pi$,
\begin{align}
D^2([P,t],[Q,s])&= t^2 + s^2 - 2st \cos (d_{PQ})\notag\\
&= (t-s)^2 + 2st (1-\cos (d_{PQ})) \notag\\
& \geq 2st(1-\cos(d_{PQ}))\label{distancelowerbnd}\\
&= st D^2(\Pi([P,t]),\Pi([Q,s])).\notag
\end{align}

\section{A monotonicity formula}\label{monotoness}

The goal of this section is to prove a proposition analogous to \cite[Proposition 3.1]{daskal-meseCAG}. The reader would benefit from familiarity with Section 3, up through Lemma 3.5, of that paper.

 Let ${\bf B}$ be a local model. In each wedge $W$, we use Euclidean coordinates $(x^1,...,x^n)$.
  For $x,y \in {\bf B}$, denote the induced Euclidean
  distance by $|x-y|$. Thus, if $x=(x^1,...,x^n)$ and
  $y=(y^1,...,y^n)$ are on the same wedge of ${\bf B}$, then
  $|x-y|^2=\sum_{i=1}^n{(x^i-y^i)^2}$.
  Let $(r,\theta_1,...,\theta_{n-1})$ denote polar coordinates, so $r$
represents radial distance from the origin and $\theta=(\theta_1,...,\theta_{n-1})$ are the standard coordinates on the $(n-1)$-sphere.

Presume, unless otherwise stated, that $g$ is a normalized Lipschitz metric defined on ${\bf B}(r)$. For $\sigma \in (0,r)$, set
\begin{equation} \label{wrtE1}
 E_g^f(\sigma)= \int_{{\bf B}(\sigma)} |\nabla f|_g^2 d\mu_g
\end{equation}
and
\begin{equation} \label{wrtE2}
I_g^f(\sigma, Q) = \int_{\partial  {\bf B}(\sigma)} d^2 (f,Q)
d\Sigma_g
\end{equation}
for  $Q \in Y$.  Here $d\Sigma_g$ is the measure on $\partial
{\bf B}(\sigma)$ induced by  $g$.  
\begin{notation}For simplicity, in the {\it{rest of this section}} we
will use the notation
\[ 
E(\sigma)= E_g^f(\sigma) \ \mbox{ and } \ I(\sigma)=I(\sigma, Q)= I_g^f(\sigma, Q),
\]
if $Q$ is a generic point. 
Furthermore in all statements we assume that the metric $g$ is normalized.
\end{notation}

We begin with a technical lemma which provides a unique center of mass for energy minimizers into sufficiently small balls. See \cite[Lemma 2.5.1]{korevaar-schoen1} for the analogous statement for $L^2$ maps with NPC targets.

\begin{lemma}\label{cenm}Let $(Y,d)$ be a CAT(1) space and $0 < \tau < \frac \pi 4$. If $f:({\bf B}(r),g) \to \B_{\tau}(P) \subset Y$ is an $L^2$ map,  then for each $0<\sigma<r$ there exists a unique $Q_\sigma \in Y$ such that
\begin{equation*}
I(\sigma, Q_\sigma)= \inf_{Q \in Y}I(\sigma,Q).
\end{equation*}
\end{lemma}
\begin{proof}Note that it is enough to consider points $Q$ in $\B_\tau(P)$ since the projection function is distance decreasing on balls of radius $\frac \pi 4$ in CAT(1) spaces. By Lemma \ref{lemmacvx}, for $x \in {\bf B}(r)$ and $Q,R \in \B_\tau(P)$, \eqref{fcvx} implies that
\[
\frac 1{8}\cos \left(d_{Q_{\frac 12}f(x)}\right) d^2_{QR} \leq \frac 12 (d^2(R,f(x)) + d^2(Q,f(x)))-d^2(Q_{\frac 12},f(x))
\]where $Q_{\frac 12}$ is the midpoint between $Q,R$. Note that $\cos \left(d_{Q_{\frac 12}f(x)}\right) \geq \cos(2\tau) >0$. Thus, integrating over $({\bf B}(\sigma),g)$ implies that 
\[
d^2_{QR} \leq \frac C{\cos(2\tau)} \left(\frac 12 \left( I(\sigma,R) + I(\sigma,Q)\right) - I(\sigma, Q_{\frac 12})\right).
\]It follows that any minimizing sequence for $I(\sigma,Q)$ is Cauchy and therefore there is a unique minimum.
\end{proof}

We next prove a type of subharmonicity result for the $d^2$ function. See \cite[Proposition 2.2.]{gromov-schoen}, \cite[Lemma 3.3]{daskal-meseCAG} for a similar result when $Y$ is NPC. Note that in the NPC setting the integral of $d^2(f, Q)|\nabla f|_g^2$ does not appear in \eqref{predomvareq}. 

 \begin{lemma} \label{predomainvariation}Let $0< \tau < \frac \pi 2$ and $f:({\bf B}(r),g) \rightarrow \B_{\tau}(P) \subset Y$ be an energy minimizing map and $(Y,d)$ a CAT(1) space. Presume that $Q \in \B_{\tau}(P)$. 
Then
 for all $0<\sigma \leq r$
 \begin{equation}\label{predomvareq}
2E(\sigma) -  \int_{{\bf B}({\sigma})}  d^2(f, Q)|\nabla f|_g^2d\mu_g\leq \int_{\partial{\bf B}({\sigma})} \langle\nabla |x|, \nabla d^2(f, Q)\rangle_g d\Sigma_g.  
\end{equation}
\end{lemma}

\begin{proof}Define 
$f_{\eta}:({\bf B}(r),g) \rightarrow Y$ 
by setting
\[
f_{\eta}(x)=(1-\eta(x)) f(x)+\eta(x) Q
\]
for $\eta \in C^{\infty}_c ({\bf B}(r))$. 
Letting $S = f(x),P = f(y),  \eta' =\eta(y)$, we use the estimate of Lemma \ref{tri1} to observe that  for $\hat d(x):= d(Q, f(x))$, 
\begin{align*}
d^2(f_\eta(y), f_\eta(x)) &\leq (1- 2\eta(x) + \eta(x)\hat d^2(x))d^2(f(x),f(y)) \\& \quad
-2(\eta(x)-\eta(y))(\hat d(x)-\hat d(y))\hat d(x)\\
& \quad +(\eta(y)-\eta(x))^2\hat d^2(x) + \eta^2(x)\mathrm{Quad}(d(f(x),f(y)), \hat d(x)-\hat d(y)) \\
&\quad + \mathrm{Cub}\left( d(f(x),f(y)), \hat d(x)-\hat d(y), \eta(x)-\eta(y) \right).
\end{align*}
Divide by $\epsilon^{n+1}$ and fix $x\in {\bf B}(r)_\epsilon$ where 
\[
{\bf B}(r)_\epsilon= \{x \in {\bf B}(r): d(x, \partial {\bf B}(r))> \epsilon\}.
\] Let $S(x,\epsilon)$ denote the $\epsilon$-sphere centered at $x$. By integrating over all $y \in S(x,\epsilon)$ with respect to the induced measure on $S(x,\epsilon)$, integrating over all $x \in {\bf B}(r)_\epsilon$, and letting $\epsilon \rightarrow 0$,
 we obtain
 \begin{eqnarray*}
 \int_{{\bf B}(r)}|\nabla f_{\eta}|_g^2 d\mu_g &  \leq &  \int_{{\bf B}(r)} |\nabla f|_g^2 d\mu_g -2\int_{{\bf B}(r)} \eta |\nabla f|_g^2 d\mu_g+ \int_{{\bf B}(r)} \eta d^2(f, Q)|\nabla f|_g^2 d\mu_g
 \\
 & &   - \int_{{\bf B}(r)} \langle \nabla \eta , \nabla d^2(f, Q)\rangle_g  d\mu_g+ O(\eta^2, |\nabla \eta|_g^2).
 \end{eqnarray*} Note that the cubic error terms either vanish as $\epsilon \to 0$ or can be absorbed into the remaining error.
 
 Now note that the energy of $f$ is bounded from above by the energy of $f_{\eta}$.  Thus, 
\begin{equation*}
2\int_{{\bf B}({r})} \eta |\nabla f|_g^2d\mu_g -  \int_{{\bf B}({r})} \eta d^2(f, Q)|\nabla f|_g^2d\mu_g\leq -\int_{{\bf B}({r})}\langle \nabla \eta , \nabla d^2(f, Q)\rangle_g d\mu_g   + O(\eta^2, |\nabla \eta|_g^2).
\end{equation*}
Replace $\eta$ by $\alpha \eta$, divide by $\alpha$ and let $\alpha \rightarrow 0$ to cancel out  the $O(\eta^2, |\nabla \eta|_g^2)$ term.  Letting $\eta$ approximate the characteristic function on ${\bf B}(\sigma)$ implies \eqref{predomvareq}. 
 \end{proof}
  
  \begin{lemma}\label{tau3lemma}Let $0< \tau < 1$ and $f:({\bf B}(r),g) \rightarrow \B_{\tau}(P) \subset Y$ be an energy minimizing map, $(Y,d)$ a CAT(1) space, and $g$ a normalized Lipschitz metric. Then, for all $0<\sigma<r$,
 \begin{equation}\label{EIinequ0}
\frac 12 E(\sigma) \leq  I(\sigma)^{1/2} \left(\left(\int_{\partial{\bf B}(\sigma)} \left|\frac{\partial f}{\partial r}\right|^2 d\Sigma_g\right)^{1/2} + c \sigma(E'(\sigma))^{1/2} \right),
\end{equation}where $c$ depends on ${\bf B}(r)$, and the Lipschitz bound and ellipticity constant of $g$. 
\end{lemma}
\begin{proof}
By \eqref{predomvareq} and the Lipschitz bound $|g^{ij}- \delta^{ij}| \leq c\sigma$, for $Q \in \B_\tau(P)$,
\begin{align}\label{UseAgain}
(2-4\tau^2)E(\sigma)&\le \int_{\partial{\bf B}({\sigma})} \langle\nabla |x|, \nabla d^2(f, Q)\rangle_g d\Sigma_g \notag \\
&=\int_{\partial{\bf B}({\sigma})} g^{ij}\frac{x_j}{|x|}\,\frac{\partial}{\partial x^i}d^2(f,{Q})\,d\Sigma_g \notag \\
&\le \int_{\partial{\bf B}({\sigma})}\frac{\partial }{\partial r}d^2(f, {Q})\, d\Sigma_g+c\sigma\int_{\partial{\bf B}({\sigma})}\sum\limits_i\Big|\frac{\partial}{\partial x^i}d^2(f, {Q})\Big|d\Sigma_g\\
& \notag = 2 \int_{\partial {\bf B}({\sigma})}d(f,Q) \frac \partial{\partial r} d(f,Q) d\Sigma_g + 2c\sigma\int_{\partial {\bf B}({\sigma})}d(f,Q) \sum_i\left|\frac\partial{\partial x^i} d(f,Q)\right| d\Sigma_g\\
& \leq  2I(\sigma)^{1/2} \left(\left(\int_{\partial{\bf B}(\sigma)} \left|\frac{\partial f}{\partial r}\right|^2 d\Sigma_g\right)^{1/2} + c \sigma(E'(\sigma))^{1/2} \right).\notag
\end{align}
In the final inequality we use H\"older's inequality and the fact that
\[
\left|\frac \partial{\partial r} d(f,Q)\right|^2 \leq \left|\frac{\partial f}{\partial r}\right|^2, \quad
 \left|\frac \partial {\partial x^i}d(f,Q)\right|^2 \leq \left|\frac{\partial f}{\partial x^i}\right|^2, \quad |\nabla f|^2 \leq \frac 1{\lambda^2}|\nabla f|^2_g,
\]where $\lambda$ is the ellipticity constant of $g$.
\end{proof}

\begin{lemma} \label{chenenergy}
Let $0<\tau<\frac \pi 4$ and $f:({\bf B}(r),g) \rightarrow \B_\tau(P)\subset Y$ be an energy minimizing map into a CAT(1) space $(Y,d)$. There exist $\sigma_0>0$ and $\gamma>0$ depending on ${\bf B}(r)$, and the Lipschitz bound and the ellipticity constant of $g$ so that 
\[
\sigma \mapsto \frac{E(\sigma)}{\sigma^{n-2+2\gamma}}, \ \sigma \in (0,\sigma_0)
\]
is non-decreasing.
\end{lemma}
\begin{proof}
 Let $Q_\sigma \in Y$ such that
\[
I(\sigma,Q_\sigma) = \inf_{Q \in Y} I(\sigma,Q)
\]where the existence is guaranteed by Lemma \ref{cenm}.
We now follow the exact argument of \cite[Lemma 3.5]{daskal-meseCAG}. Note that their invocation of \cite[(3.12)]{daskal-meseCAG} is replaced by \eqref{EIinequ0} here. All other inequalities they reference arise from appropriate domain variations and are therefore true for maps into $Y$.
\end{proof}

As in \cite{daskal-meseCAG}, all of the previous results extend from the setting of normalized metrics to admissible complexes with Lipschitz Riemannian metrics. See \cite[p. 289-290]{daskal-meseCAG} to understand how the properties of the map $L_x$ in \cite[Proposition 2.1]{daskal-meseCAG} affect the energy of $f$ and the domain over which Lemma \ref{chenenergy} can be applied.


Let ${\bf B}$ be a dimension-$n$, codimension-$(n-k)$ local model and $g$  a  Lipschitz metric on ${\bf B}(r)$ with ellipticity
constant $\lambda \in (0,1]$.  For  $x \in {\bf B}(r)$, 
let $R(x)$
denote the radius of the largest homogeneous ball centered at $x$ contained in  ${\bf B}(r)$. 
 The value  $\sigma_0>0$  was defined above as the upper bound for which the monotonicity  formula of Lemma~\ref{chenenergy} holds for any energy minimizing map from a  local model with a normalized metric.  Therefore, the monotonicity formula  for $f \circ L_x$ is valid for balls ${\bf B}'(\sigma)$ contained in ${\bf B}'(r_0(x))$ where 
 \begin{equation} \label{defr_0}
 r_0(x):=\min\{\sigma_0, \lambda R(x)\}.
 \end{equation}    
Recalling that $B_x(\sigma)$ is the Euclidean ball about $x$ of radius $\sigma$, we define  $E_x(\sigma)$ for $\sigma$
sufficiently small by setting
\[
E_x(\sigma)=\int_{B_x(\sigma)} |\nabla f |^2
d\mu_g.
\]

\begin{proposition} \label{relax}
Let ${\bf B}$ be a dimension-$n$, codimension-$\nu$ local model, $g$ a  Lipschitz Riemannian metric  defined on ${\bf B}(r)$ with ellipticity constant $\lambda \in (0,1]$, $(Y,d)$ a CAT(1) space and $f:({\bf B}(r),g) \rightarrow \B_\tau(P) \subset Y$  an energy minimizing map. If $0<\tau< \frac \pi 4$, then there exist  constants $\gamma>0$ and $C \geq 1$ depending on ${\bf B}(r)$, the  Lipschitz bound and the ellipticity constant of $g$ so that for every $x \in {\bf B}(r)$, 
\begin{equation} \label{gooddecay1}
\frac{E_x(\sigma)}{\sigma^{n-2+2\gamma}} \leq C \frac{E_x(\rho)}{ \rho^{n-2+2\gamma}}, \ 0<\sigma<\rho \leq  r(x)
\end{equation}
where
\begin{equation} \label{defr(x)}
r(x):=\lambda r_0(x)=\min\{\lambda \sigma_0, \lambda^2 R(x)\}.
\end{equation}
Here, $R(x)$ is defined as above 
and $\sigma_0>0$ is as in Lemma \ref{chenenergy}.
\end{proposition}

\begin{proof} The proof 
proceeds exactly as in the proof of \cite[Proposition 3.1]{daskal-meseCAG}, using Lemma \ref{chenenergy}.
\end{proof}

\section{H\"older Regularity}\label{Holderss}

The goal of this section is to prove Theorem \ref{Holder}. The proof is modeled on the proof of H\"older regularity in \cite{daskal-meseCAG} for minimizing maps into an NPC space. 
The reader would benefit from a familiarity with 
Section 4 of that paper. The method of proof is classical, as the regularity result will follow from a Campanato theorem and the monotonicity given by \eqref{gooddecay1}. Many of the technical aspects of this argument in
\cite{daskal-meseCAG} are related to the singular nature of the domain and thus can be immediately applied for CAT(1) targets.
We will highlight the key places where the target curvature plays a role and provide suitable adaptations of the arguments involved.

Before proceeding to the main argument, we prove a technical lemma.
\begin{lemma}\label{L2Close}Let $\Omega$ be a Euclidean domain, $(Y,d)$ a CAT(1) space, and $0<\tau<\frac \pi 4$. If $f:\Omega\to \B_{\tau}(P) \subset Y$ is an $L^2$ map  then, for all $\varepsilon>0$, there exists $h_\epsilon:\Omega_\epsilon \to Y$ Lipschitz such that 
\[
\int_{\Omega_\epsilon}d^2(f,h_\epsilon) d\mu_g< \epsilon,
\]where $\Omega_\epsilon:= \{x \in \Omega: d(x,\partial \Omega) >\epsilon\}$.
\end{lemma}
\begin{proof}For $f:\Omega \rightarrow Y$, recall  $\overline{f}:\Omega \rightarrow {\mathcal C}Y$ is defined by setting
 \[
\overline{f}(x)=[f(x),1].
\](See Section \ref{conesection} for further relevant definitions and energy comparisons.)

We  use the notation $\B^{{\mathcal C}Y}_r(\cdot)$ to denote a ball of radius $r$ in ${\mathcal C}Y$.
By \eqref{D_d_compare},
\[
\overline{f}(x) \in \B^{{\mathcal C}Y}_{\tau}([P,1]), \ \ \forall x \in \Omega .
\]and  $\overline{f}$ is an $L^2$-map into $\mathcal CY$. Since ${\mathcal C}Y$ is NPC, we can apply the mollification procedure of \cite[Section 1.5]{korevaar-schoen2} to produce a Lipshitz map $g_{\epsilon}:  \Omega_{\epsilon} \rightarrow \B^{{\mathcal C}Y}_{\tau}([P,1]) \subset {\mathcal C}Y$ such that
\[
\int_{\Omega_{\epsilon}} D^2(\overline{f},g_{\epsilon}) dx < \frac{\epsilon}{4}.
\]
Write $g_{\epsilon}(x)=[\varphi(x),t(x)]$.  The map $g_{\epsilon}(x)$ is constructed as the center of mass of the map $\overline{f}$ with respect to a probability measure $\eta_{\epsilon}(x-y)dy$ where $\eta_{\epsilon}$ can be chosen to be a function  with compact support in a small ball centered at $0$. Therefore, since Image$(\overline{f}) \subset Y \times \{1\}$, we can assume that $D(g_{\epsilon}(x),Y \times \{1\})$ satisfies 
\[
|1-t(x)| < \left(\mbox{volume}(\Omega_\epsilon) \right)^{-1/2}\sqrt{\frac{\epsilon}{4}}.
\] 
Thus, 
\[
\int_{\Omega_{\epsilon}} D^2(\overline{f},\Pi \circ g_{\epsilon}) dx \leq \int_{\Omega_{\epsilon}} D^2(\overline{f},g_{\epsilon}) dx  + \int_{\Omega_{\epsilon}}D^2(g_{\epsilon}, \Pi \circ g_{\epsilon}) dx <\frac {\epsilon}{2}
\] 
where $\Pi:  {\mathcal C}Y \rightarrow Y \times \{1\}$ is the projection map as in Section~\ref{conesection}. Since $|1-t(x)|$ is bounded for all $x \in \Omega$ and $g_\epsilon$ is Lipschitz, \eqref{distancelowerbnd} implies that $\Pi\circ g_\epsilon$ is Lipschitz on $\Omega_\epsilon$.
Define
\[
h_{\epsilon}:  \Omega_{\epsilon} \rightarrow Y, \ \ \ \ h_{\epsilon}=\Pi \circ g_{\epsilon}
\] 
by identifying $Y$ with $Y \times \{1\} \subset {\mathcal C}Y$.  
Then by \eqref{D_d_compare}
\begin{eqnarray*}
\int_{\Omega_{\epsilon}} d^2(f,h_{\epsilon}) dx & \leq &2 \int_{\Omega_{\epsilon}}D^2(\overline{f},\Pi \circ g_{\epsilon}) dx < \epsilon.
\end{eqnarray*}
\end{proof}

We now prove a Campanato type lemma. In \cite[Lemma 4.1]{daskal-meseCAG}, the authors prove a similar result for any $L^2$ map into an NPC space.
\begin{lemma}\label{campanato}
Let ${\bf B}$ be a dimension-$n$, codimension-$\nu$ local model, $g$ a Lipschitz Riemannian metric on ${\bf B}(r)$,
$(Y,d)$ a CAT(1) space, $0< \tau < \frac \pi 4$ and  $f:({\bf B}(r),g) \rightarrow \B_\tau(P) \subset Y$  an  $L^2$ map. Fix $\varrho \in (0,1)$. If there exist $K>0$, $R \in (0,(1-\varrho)r)$ and $\beta \in (0,1]$ such that
\begin{equation} \label{energygrowth1}
\inf_{Q \in Y} \sigma^{-n}  \int_{B_x(\sigma)}d^2(f,Q ) \ d\mu_g  \leq K^2 \sigma^{2\beta}, \forall x
\in {\bf B}(\varrho r) \mbox{ and } \sigma \in (0,R),
\end{equation}
then there exists $C>0$ and a representative in the $L^2$-equivalence class of $f$, which we still denote by $f$, such that
\[
d(f(x),f(y)) \leq C|x-y|^{\beta}, \ \forall x,y \in {\bf B}(\varrho r)
\]
with $C$ depending on $K$,  $r$, $R$, $\beta$, $\varrho$ and ${\bf B}(r)$.
\end{lemma}
\begin{proof}
The lemma will follow from the Campanato lemma \cite[Lemma 4.1]{daskal-meseCAG}, provided that each aspect of the proof that relied on the non-positive curvature of the target still holds if the target is CAT(1) and $f$ has small image. The NPC hypothesis gave the existence and uniqueness of $Q_{x,\sigma}$ for each $x \in {\bf B}(\varrho r)$ and $\sigma \in (0,R)$. Lemma \ref{cenm} above provides this for our setting. The NPC condition also provided the existence of Lipschitz maps $L^2$ close to $f$. For a CAT(1) space $Y$, we appeal to Lemma \ref{L2Close} above, since, by hypothesis, $f$ has small image. All other aspects of the proof are related to properties of the domain, and thus carry through with no trouble. 
\end{proof}

Recall the following proposition \cite[Proposition 4.3]{daskal-meseCAG}, which converts the monotonicity information of \eqref{gooddecay1} into a uniform estimate on the decay of the scale invariant energy for all $x$.
\begin{proposition} \label{propjoe}
Let ${\bf B}$ be a dimension-$n$, codimension-$\nu$  local model, $g$ a Lipschitz Riemannian metric defined on ${\bf B}(r)$,  $(Y,d)$ a metric space and $f:({\bf B}(r),g) \rightarrow Y$ a finite energy map.  Fix $\varrho \in (0,1)$ and suppose that for $x \in {\bf B}\left(\varrho r \right)$ there exist $\beta>0$ and $\hat{C}\geq 1$ so that
\begin{equation} \label{roth}
\frac{E_x(\sigma)}{\sigma^{n-2+2\beta}} \leq \hat{C} \frac{E_x(\rho)}{\rho^{n-2+2\beta}} ,\  \ 0 < \sigma \leq \rho \leq r(x)
\end{equation}
where $r(x)$ is as defined in \eqref{defr(x)}.
Then
there exist $K$ and $R>0$ depending only on the total energy of $f$, $E^f$,  the ellipticity
constant and Lipschitz bound of $g$,  ${\bf
B}(r)$ and $\varrho$ so that
\[
E_x(\sigma) \leq K^2 \sigma^{n-2+2\beta}, \hspace{0.1in} \forall x
\in {\bf B}(\varrho r), \sigma<R.
\]
 \end{proposition}

This immediately implies the H\"older regularity for a local model.
\begin{theorem} \label{chen}
Let  ${\bf B}$ be a local model, $g$  a 
Lipschitz Riemannian metric defined on ${\bf B}(r)$, $(Y,d)$ a CAT(1) space and $f:({\bf B}(r),g) \rightarrow \B_\tau(P) \subset Y$  an energy minimizing map where $0<\tau< \frac \pi 4$.  For $\varrho \in(0,1)$,    there exist $C_H>0$ and
$\gamma>0$  depending only on the Lipschitz bound and ellipticity constant of $g$, $E^f$, ${\bf B}(r)$ and $\varrho$ such  that
\[
d(f(x),f(y)) \leq C_H |x-y|^{\gamma}, \hspace{0.1in} \forall x,y \in
{\bf B}(\varrho r).
\]
\end{theorem}

\begin{proof}
The result follows immediately from \eqref{gooddecay1},  Proposition~\ref{propjoe}, the Poincar\'e inequality of \cite[Theorem 2.7]{daskal-meseCAG}, and Lemma~\ref{campanato}.
\end{proof}

By using \cite[Proposition 2.1]{daskal-meseCAG} we obtain Theorem \ref{Holder}.

\section{Improved H\"older Regularity}\label{Improved}

To extend the regularity from H\"older to Lipschitz requires a better result than Theorem \ref{chen} provides. The objective of this section is two-fold. First, we prove that the order function $\mathrm{ord}^f(x)$ is well-defined (see Proposition \ref{orderdef} and Definition \ref{orderdef2}). Second, we use monotonicity to demonstrate that the H\"older regularity on a ball can be improved to have H\"older exponent equal to $\alpha\leq \mathrm{ord}^f(x)$ for $x \in {\bf B}(r)$.
The Lipschitz regularity will then immediately hold in any neighborhood with $\alpha \geq 1$. 

In \cite{daskal-meseCAG}, the authors proved the stronger H\"older regularity in parallel with the weaker version. In the CAT(1) setting, however, we rely in a fundamental way on the weaker H\"older result. We use the weak H\"older result in \eqref{predomvareq} to improve the inequality from \eqref{UseAgain}. This improvement allows us adapt the techniques of \cite{daskal-meseCAG} to our setting. Following their ideas, we demonstrate that the order function is well-defined. We then demonstrate that $\frac{E(\sigma)}{\sigma^{n-2+\alpha}}$ is monotone, which immediately implies the improved regularity.

Throughout this section, unless explicitly stated otherwise, presume that $g$ is a normalized Lipschitz metric on ${\bf B}(r)$.

\subsection{The order function}The goal of this subsection is to prove that the order $\alpha:= \lim_{\sigma\to 0^+} \frac{\sigma E(\sigma)}{I(\sigma)}$ exists. In the Euclidean setting, the existence of the limit follows from proving the differential inequality $\frac{E'(\sigma)}{E(\sigma)} - \frac{I'(\sigma)}{I(\sigma)} + \frac 1\sigma \geq 0$, which implies that the function $\sigma \mapsto \frac{\sigma E(\sigma)}{I(\sigma)}$ is monotone. Under the current hypotheses, we cannot hope to prove a differential inequality of exactly the desired type. The inequality we determine includes additional terms. Nevertheless, we still show that the limit $\alpha$ exists.

We begin by recalling two essential inequalities derived in \cite[(3.9),(3.17)]{daskal-meseCAG} for energy minimizing maps from a local model into a metric space target. These calculations use only domain variations and the Lipschitz assumption on the domain metric and thus immediately extend to our setting.

\begin{lemma}
Let $f:({\bf B}(r),g) \to Y$ be an energy minimizing map, $Y$ a metric space, and $g$ a normalized Lipschitz metric. Then there exist constants $c_1,c_2>0$ and $\sigma_0>0$ small, all depending only on ${\bf B}(r)$ and the Lipschitz bounds of $g$, such that for all  $0<\sigma\leq\sigma_0$ and $Q \in Y$,
\begin{equation*}\label{e02}
\Big(1+c\sigma\Big)\frac{E'(\sigma)}{E(\sigma)}\ge \frac{n-2}{\sigma}+\frac{2}{E(\sigma)}\int_{\partial {\bf B}(\sigma)}\Big|\frac{\partial f}{\partial r}\Big|^2\, d\Sigma_g-c_1,
\end{equation*}
and
\begin{equation}\label{e03}
\left|\frac{I'(\sigma)}{I(\sigma)}- \frac{n-1}{\sigma}-\frac{1}{I(\sigma)}\int_{\partial{\bf B}(\sigma)}\frac{\partial}{\partial r}d^2(f, Q)\, d\Sigma_g\right|\le c_2.
\end{equation}
\end{lemma}

Therefore, for $c_3=c_1+c_2$,
\begin{align}\label{e04}
\notag\Big(1+c\sigma\Big)\frac{E'(\sigma)}{E(\sigma)}&-\frac{I'(\sigma)}{I(\sigma)}+\frac{1}{\sigma}+c_3\\ &\ge \frac{2}{E(\sigma)}\int_{\partial{\bf B}(\sigma)}\Big|\frac{\partial f}{\partial r}\Big|^2\, d\Sigma_g-\frac{1}{I(\sigma)}\int_{\partial {\bf B}(\sigma)}\frac{\partial}{\partial r}d^2(f, Q)\, d\Sigma_g\\
&=\frac{1}{E(\sigma)I(\sigma)}\Big(2I(\sigma)\int_{\partial{\bf B}(\sigma)}\Big|\frac{\partial f}{\partial r}\Big|^2\, d\Sigma_g-E(\sigma)\int_{\partial {\bf B}(\sigma)}\frac{\partial}{\partial r}d^2(f, Q)\, d\Sigma_g\Big).\notag
\end{align}
We now determine a lower bound for the right hand side of \eqref{e04}, modifying the differential inequality to one more conducive to the proof of monotonicity.

 \begin{lemma} \label{domainvariationDM}
Suppose that $\bf{B}$ is a local model, $(Y,d)$ is a CAT(1) space, and $g$ is a normalized Lipschitz metric. If $f:({\bf B}(r),g) \to {\B}_{\tau}(P) \subset Y$ is an energy minimizing map with $0<\tau< \frac \pi 4$, then
 for ${\bf B}({\sigma}) \subset {\bf B}(r/2)$ and any $Q \in \overline{\B_{C_H \sigma^\gamma}(f(0))}$, 
\begin{equation}\label{EIinequ}
2\Big(1-c'\sigma^{2\gamma}\Big)E(\sigma) \leq  \int_{\partial {\bf B}(\sigma)} \frac \partial{\partial r} d^2(f,Q) d\Sigma_g + I(\sigma) +k\sigma^2E'(\sigma)
\end{equation}
and
\begin{equation}\label{monostart}
E(\sigma)\int_{\partial{\bf B}({\sigma})} \frac{\partial }{\partial r}d^2(f, Q)d\Sigma_g\le 2I(\sigma)\int_{\partial{\bf B}({\sigma})}\Big|\frac{\partial f}{\partial r}\Big|^2d\Sigma_g+\big(c\sigma+c c'\sigma^{2\gamma}\big)I(\sigma)E'(\sigma)
\end{equation}
where $c,k$ depend on ${\bf B}(r)$, the Lipschitz bound and the ellipticity constant of $g$ and $c', \gamma$ depend on $E^f$, ${\bf B}(r)$, and the Lipschitz bound and ellipticity constant of $g$.
\end{lemma}
\begin{proof}
First observe that by Theorem \ref{chen}, $f(\partial {\bf B}(\sigma)) \subset \overline{\B_{C_H \sigma^\gamma}(f(0))}$. By \eqref{predomvareq} and the Lipschitz bound $|g^{ij}(x)-\delta^{ij}| \leq c\sigma$, for $|x| \leq \sigma$, we can improve the estimate in \eqref{UseAgain} to
\begin{align}\label{line1}
\Big(2-C_H^2\sigma^{2\gamma}\Big)E(\sigma)
&\le \int_{\partial{\bf B}({\sigma})}\frac{\partial }{\partial r}d^2(f, Q)\, d\Sigma_g+c\sigma\int_{\partial{\bf B}({\sigma})}\sum\limits_i\Big|\frac{\partial}{\partial x_i}d^2(f, Q)\Big|d\Sigma_g
\end{align}
where $C_H$ is the H\"older constant, and $\gamma$ is the H\"older exponent.

Then, for $c'=C_H^2/2$, 
\eqref{EIinequ} follows by applying the following elementary inequality to the last term in \eqref{line1}
\begin{align*}
2c\sigma d(f, Q)\sum_{i=1}^n\Big|\frac{\partial}{\partial x_i}d(f, Q)\Big| &\leq  d^2(f,Q) + c^2\sigma^2 \sum_{i=1}^n\Big|\frac{\partial}{\partial x_i}d(f, Q)\Big|^2
\\ & \leq d^2(f,Q) + c^2\sigma^2 \sum_{i=1}^n\Big|\frac{\partial f}{\partial x_i}\Big|^2.
\end{align*}
To prove \eqref{monostart}, first note that if $\int_{\partial{\bf B}({\sigma})} \frac{\partial }{\partial r}d^2(f, Q)d\Sigma_g \leq 0$, the result holds simply because the right hand side of the inequality is non-negative. So suppose that $\int_{\partial{\bf B}({\sigma})} \frac{\partial }{\partial r}d^2(f, Q)d\Sigma_g \geq 0$.
Recall the estimates determined in \eqref{UseAgain}:
\begin{equation}\label{e05a}
\int_{\partial{\bf B}({\sigma})} \frac{\partial }{\partial r}d^2(f, Q)d\Sigma_g \leq 2 \left(I(\sigma)\int_{\partial{\bf B}({\sigma})}\Big|\frac{\partial f}{\partial r}\Big|^2d\Sigma_g\right)^{1/2} \leq 2cI(\sigma)^{1/2}E'(\sigma)^{1/2},
\end{equation}
\begin{equation}\label{e05b}
 \int_{\partial{\bf B}({\sigma})} \sum\limits_i\Big|\frac{\partial }{\partial x_i}d^2(f,Q)\Big|d\Sigma_g \leq 2cI(\sigma)^{1/2}E'(\sigma)^{1/2}.
\end{equation}Note that in the above equations, $c$ depends on the ellipticity constant of $g$. In what follows, $c$ may increase from one line to the next, but its dependence will always be only on ${\bf B}(r)$, the Lipschitz bound and the ellipticity constant of $g$.

Using \eqref{line1}, \eqref{e05a}, \eqref{e05b} we observe that
\begin{equation}\label{e06}
\begin{split}
2\big(1-c'&\sigma^{2\gamma}\big)E(\sigma)\int_{\partial{\bf B}({\sigma})} \frac{\partial }{\partial r}d^2(f, Q)d\Sigma_g\\
&\le \Big(\int_{\partial{\bf B}({\sigma})} \frac{\partial }{\partial r}d^2(f, Q)d\Sigma_g
 +c\sigma \int_{\partial{\bf B}({\sigma})} \sum\limits_i\Big|\frac{\partial }{\partial x_i}d^2(f,Q)\Big|d\Sigma_g\Big)\,\int_{\partial{\bf B}({\sigma})}  \frac{\partial }{\partial r}d^2(f, Q)d\Sigma_g\\
& \le c(1+\sigma)I(\sigma)E'(\sigma).
\end{split}
\end{equation}
Thus, for sufficiently small $\sigma>0$,
$$
E(\sigma)\int_{\partial{\bf B}({\sigma})}\frac{\partial }{\partial r}d^2(f, Q)d\Sigma_g\le c I(\sigma)E'(\sigma).
$$ 
Now, using the middle inequality in \eqref{e05a} and substituting the above inequality into \eqref{e06} implies that
\begin{align*}
E(\sigma)\int_{\partial{\bf B}({\sigma})} &\frac{\partial }{\partial r}d^2(f, Q)d\Sigma_g\\
&\leq 2I(\sigma) \int_{\partial{\bf B}({\sigma})}\Big|\frac{\partial f}{\partial r}\Big|^2d\Sigma_g+ c\sigma I(\sigma) E'(\sigma) + 2c'\sigma^{2\gamma}E(\sigma)\int_{\partial{\bf B}({\sigma})} \frac{\partial }{\partial r}d^2(f, Q)d\Sigma_g\\
& \leq  2I(\sigma) \int_{\partial{\bf B}({\sigma})}\Big|\frac{\partial f}{\partial r}\Big|^2d\Sigma_g+ c \sigma I(\sigma) E'(\sigma) +
cc'\sigma^{2\gamma} I(\sigma)E'(\sigma).
\end{align*}
\end{proof}

Combining \eqref{e04} and \eqref{monostart}, we conclude that for sufficiently small $\sigma$, 
\begin{equation*}
\Big(1+c\sigma+c c'\sigma^{2\gamma}\Big)\frac{E'(\sigma)}{E(\sigma)}-\frac{I'(\sigma)}{I(\sigma)}+\frac{1}{\sigma}+c_3\ge 0.
\end{equation*}
Note that if $\gamma\geq \frac 12$, we may appeal directly to the work of \cite{daskal-meseCAG} since the term $c\sigma$ dominates. Therefore, we presume that $\gamma<1/2$. In what follows, for notational simplicity, we rescale the domain metric $g$ so that 
$c \le 1$, since $c$ depends only on the domain metric. If we assume that $\sigma_0=1$ and let $C\geq1+2c'$ then
\begin{equation}\label{e08}
(1+C\sigma^{2\gamma})\frac{E'(\sigma)}{E(\sigma)}-\frac{I'(\sigma)}{I(\sigma)}+\frac{1}{\sigma}+c_3\ge 0.
\end{equation}
For the analogous inequality in the NPC setting see \cite[(3.20)]{daskal-meseCAG}.

We remark that due to the extra term $\sigma^{2\gamma}$, the original monotonicity in \cite{gromov-schoen} no longer works. Following the ideas of \cite{daskal-meseCAG}, we introduce a modified energy. Let 
$$
J(\sigma)=\max\limits_{s\in(0, \sigma]}I(s),
$$
and 
$$
A:=\Big\{\sigma:\frac{E'(\sigma)}{E(\sigma)}-\frac{J'(\sigma)}{J(\sigma)}+\frac{1}{\sigma}+c_3\le 0\Big\}.
$$Note that $A$ is exactly the set on which a standard monotonicity formula fails. 
For $\sigma\in(0, 1)$, we define the modified energy
\begin{equation}\label{Fdef}
F(\sigma)=E(\sigma)\exp\Big(\varphi(\sigma)\Big),
\end{equation}
where
$$
\varphi(\sigma)=-\int_{A\cap(\sigma, 1)}Cs^{2\gamma}\frac{E'(s)}{E(s)}ds
$$and $C$ is as in \eqref{e08}. 
Then exactly as in the proof of \cite[Lemma 3.7]{daskal-meseCAG},
\begin{equation}\label{diffF}
\frac{F'(\sigma)}{F(\sigma)}=\left\{
\begin{array}{ll}
\frac{E'(\sigma)}{E(\sigma)}&\text{ if }\sigma\notin A\\
(1+C\sigma^{2\gamma})\frac{E'(\sigma)}{E(\sigma)}&\text{ if }\sigma\in A.
\end{array}\right.
\end{equation}
By \eqref{e08} and \eqref{diffF}, we observe that 
\begin{equation}\label{Fnondec}
\sigma\mapsto e^{c_3\sigma}\frac{\sigma F(\sigma)}{I(\sigma)}
\end{equation}
is monotone nondecreasing for any $Q \in \overline{\B_{C_H \sigma^\gamma}(f(0))}$. For $\sigma>0$ sufficiently small, $C_H \sigma^\gamma < \tau< \frac \pi 4$ and thus the projection map $\pi_\sigma: \B_\tau(P) \to \overline{\B_{C_H \sigma^\gamma}(f(0))}$ is distance decreasing. It follows that for every $\sigma>0$ sufficiently small, $Q_\sigma \in \overline{\B_{C_H \sigma^\gamma}(f(0))}$. Thus applying \eqref{Fnondec} for $\sigma_1 < \sigma_2$ sufficiently small, and noting that
by the definition of $I(\sigma, Q_\sigma)$,  $I(\sigma_2, Q_{\sigma_1})\ge I(\sigma_2, Q_{\sigma_2})$, we observe that
$$
e^{c_3\sigma_1}\frac{\sigma_1 F(\sigma_1)}{I(\sigma_1, Q_{\sigma_1})}\le e^{c_3\sigma_2}\frac{\sigma_2 F(\sigma_2)}{I(\sigma_2, Q_{\sigma_1})}\le e^{c_3\sigma_2}\frac{\sigma_2 F(\sigma_2)}{I(\sigma_2, Q_{\sigma_2})}.
$$
Therefore, $\sigma\mapsto e^{c_3\sigma}\frac{\sigma F(\sigma)}{I(\sigma, Q_\sigma)}$ is monotone and
$\lim\limits_{\sigma\to 0^+} e^{c_3\sigma}\frac{\sigma F(\sigma)}{I(\sigma, Q_\sigma)}$ exists. To show that $\lim_{\sigma \to 0^+}e^{c_3\sigma} \frac{\sigma E(\sigma)}{I(\sigma,Q_\sigma)}$ exists, it is therefore enough to consider
$$
\lim\limits_{\sigma\to 0^+}\frac{E(\sigma)}{F(\sigma)}=\lim\limits_{\sigma\to 0^+}\exp\big(-{\varphi(\sigma)}\big).
$$

\begin{lemma}\label{estbadset}
$\lim\limits_{\sigma\to 0^+}-\varphi(\sigma)=\lim\limits_{\sigma\to 0^+}\int_{A\cap (\sigma, 1)}Cs^{2\gamma}\frac{E'(s)}{E(s)}ds<\infty.$
\end{lemma}
\begin{proof}
The proof follows from straightforward modifications of the argument in \cite[Lemma  3.8]{daskal-meseCAG}.

By the definition of $A$, for all $s\in A$
\begin{equation}\label{e010}
\frac{E'(s)}{E(s)}\le \frac{J'(s)}{J(s)}.
\end{equation}
So it suffices to show that
$$
\lim\limits_{\sigma\to 0^+}\int_{A\cap (\sigma, 1)}Cs^{2\gamma}\frac{J'(s)}{J(s)}ds<\infty.
$$
Following \cite[Proof of Lemma 3.8]{daskal-meseCAG}, there exists a sufficiently large constant $C'$ depending on the domain such that for any $\epsilon>0$ and  $0<\theta_1<\theta_2\le 1$,
\begin{equation}\label{e10}
\Big[1-\big(\frac{1}{\epsilon}+\frac{C'}{\theta_1}\big)\big(\theta_2-\theta_1\big)\Big]J(\theta_2)-\epsilon ME(\theta_2)\le J(\theta_1)
\end{equation}
where 
$$
 \frac 1{E(\sigma)} \int_{\bf{B}(\sigma)}|\nabla f|^2d\mu_g \leq M.
$$

Let 
$$
\phi(\theta, n, j)=\frac{1}{2}-\Big(2K\theta^{-jp(\theta, n)}+\frac{C'}{\theta}\Big)\Big(1-\theta\Big),
$$
where 
$$
p(\theta, n)=C\frac{1-\theta^{2\gamma n}}{1-\theta^{2\gamma}},  \quad \quad K=M e^{c_3} \frac{F(1)}{J(1)}.
$$
Note that $p(\theta, 0)=0$ and
$$
\lim\limits_{\theta\to 1^-}\phi(\theta, n, j)=\frac{1}{2}\, \,\,\text{uniformly in $j$, $n$.}
$$
Therefore, there exists $\theta_0<1$ sufficiently close to 1, such that $\phi(\theta_0, n, j)>\frac{1}{4}$.
We also choose $j$ such that $\theta_0^j<\frac{1}{4}$. Then for all $n$, we have that
\begin{equation}\label{e11}
\theta_0^j<\phi(\theta_0, n, j)=\frac{1}{2}-\Big(2K\theta_0^{-jp(\theta_0, n)}+\frac{C'}{\theta_0}\Big)\Big(1-\theta_0\Big).
\end{equation}
Now, for the chosen $\theta_0$ and $j$ as above, we let
$$
\epsilon=\frac{1}{2K}\theta_0^{jp(\theta_0, n)+n}.
$$
Then by \eqref{e10}, 
\begin{equation}\label{e12}
\Big[1-\big(2K\theta_0^{-jp(\theta_0, n)-n}+\frac{C'}{\theta_1}\big)\big(\theta_2-\theta_1\big)\Big]J(\theta_2)-\frac{E(\theta_2)J(1)}{2e^{c_3}F(1)}\theta_0^{jp(\theta_0, n)+n}\le J(\theta_1).
\end{equation}

\begin{claim}
For any $n$, 
\begin{equation}\label{Jest}
\theta_0^jJ(\theta_0^n)<J(\theta_0^{n+1}),
\end{equation}
and
\begin{equation}
\int_{\theta_0^n}^1Cs^{2\gamma}\frac{J'(s)}{J(s)} ds\le \log \theta_0^{-jp(\theta_0, n)}\leq C(\theta_0,\gamma,j).
\end{equation}
\end{claim}
Note that the proof of the lemma will follow once we prove the claim since 
\[
\lim_{\sigma \to 0^+} \int_{A \cap (\sigma,1)}Cs^{2\gamma}\frac{E'(s)}{E(s)} ds\leq \lim_{\sigma \to 0^+} \int_{A \cap (\sigma,1)}Cs^{2\gamma}\frac{J'(s)}{J(s)} ds\leq \lim_{n \to \infty} \int_{\theta_0^n}^1Cs^{2\gamma}\frac{J'(s)}{J(s)} ds
\]and by the claim, the right hand side is bounded independent of $n$.
\end{proof}
\begin{proof}[Proof of claim]

We proceed by induction on the powers of $\theta_0$.
First, take $n=0$, $\theta_1=\theta_0$, $\theta_2=1$, and notice that $F(1)=E(1)$ by the definition of $F(\sigma)$, and thus $\frac{E(1)J(1)}{e^{c_3}F(1)}=\frac{J(1)}{e^{c_3}}\le J(1)$. Then by \eqref{e11} and \eqref{e12},
$$
\theta_0^jJ(1)<\phi(\theta_0, 0, j)J(1)=\Big[\frac{1}{2}-\big(2K+\frac{C'}{\theta_0}\big)(1-\theta_0)\Big]J(1)\le J(\theta_0).
$$
Next, we assume $\theta_0^jJ(\theta_0^k)<J(\theta_0^{k+1})$ for all $k=0, 1, 2, \dots, n-1$.

By the definition of $F(\sigma)$ and \eqref{e010}, 
\begin{equation}\label{e016}
\begin{split}
\log\frac{E(\theta_0^n)}{F(\theta_0^n)}=-\varphi(\theta_0^n)&=\int_{A\cap(\theta_0^n, 1)}Cs^{2\gamma}\frac{E'(s)}{E(s)}ds\\
&\le\int_{\theta_0^n}^1Cs^{2\gamma}\frac{J'(s)}{J(s)}ds.
\end{split}
\end{equation}
We estimate
\begin{equation*}\label{IntegralJEst}
\begin{split}
\int_{\theta_0^n}^1Cs^{2\gamma}\frac{J'(s)}{J(s)}ds &= \sum\limits_{k=0}^{n-1}\int_{\theta_0^{k+1}}^{\theta_0^k}Cs^{2\gamma}\frac{d}{ds}\log J(s)ds\\
&\le \sum\limits_{k=0}^{n-1} C\theta_0^{2\gamma k}\log\frac{J(\theta_0^k)}{J(\theta_0^{k+1})}\\
&\le \sum\limits_{k=0}^{n-1} C\theta_0^{2\gamma k}\log \theta_0^{-j}\hspace{.2in}\text{(by the induction hypothesis)}\\
&=\log \theta_0^{-jC\frac{1-\theta_0^{2\gamma n}}{1-\theta_0^{2\gamma}}}=\log\theta_0^{-jp(\theta_0, n)}.
\end{split}
\end{equation*}
So since \begin{equation*}\label{e011}
\frac{\sigma F(\sigma)}{J(\sigma)} \leq e^{c_3\sigma}\frac{\sigma F(\sigma)}{J(\sigma)}\leq e^{c_3}\frac{F(1)}{J(1)}
\end{equation*}
\eqref{e016} and the integral estimate imply that
$$
\theta_0^{jp(\theta_0,n)}E(\theta_0^n)\le F(\theta_0^n)\le \frac{e^{c_3} F(1)J(\theta_0^n)}{J(1)\theta_0^n}
$$
That is
\[
\theta_0^{jp(\theta_0,n)+n}\frac{E(\theta_0^n)J(1)}{e^{c_3}F(1)} \leq J(\theta_0^n).
\]
We now take $\theta_1=\theta_0^{n+1}$, $\theta_2=\theta_0^n$ in \eqref{e12}, and together with \eqref{e11} and the above inequality we conclude that
$$
\theta_0^jJ(\theta_0^n)<\phi(\theta_0, n, j)J(\theta_0^n)=\Big[\frac{1}{2}-\Big(2K\theta_0^{-jp(\theta_0, n)}+\frac{C'}{\theta_0}\Big)\Big(1-\theta_0\Big)\Big]J(\theta_0^n)\le J(\theta_0^{n+1}).
$$
This implies that \eqref{Jest} is true for all $n$. Therefore, we may make the substitution in the integral estimate to conclude that for all $n$ 
\begin{equation}\label{e019}
\int_{\theta_0^n}^1Cs^{2\gamma}\frac{J'(s)}{J(s)}ds \le \log \theta_0^{-jC\frac{1-\theta_0^{2\gamma n}}{1-\theta_0^{2\gamma}}}=\log\theta_0^{-jp(\theta_0, n)}.
\end{equation}
Since $p(\theta_0,n)$ is increasing in $n$ and $\lim_{n \to \infty}p(\theta_0, n)=C\frac{1}{1-\theta_0^{2\gamma}}=C(\theta_0,\gamma)$, we prove the integral estimate.
\end{proof}

\begin{proposition}\label{orderdef}
Suppose that $\bf{B}$ is a local model, $(Y,d)$ is a CAT(1) space, and $g$ is a normalized Lipschitz metric. If $f:({\bf B}(r),g) \to {\B}_{\tau}(P) \subset Y$ is an energy minimizing map where $0<\tau< \frac \pi 4$, then the order 
$$
\alpha:=\lim\limits_{\sigma\to 0}\frac{\sigma E(\sigma)}{I(\sigma, Q_\sigma)} <\infty
$$
is well defined.
\end{proposition}
\begin{proof}

The monotonicity of $\sigma\mapsto \frac{\sigma F(\sigma)}{I(\sigma, Q_\sigma)}$ together with Lemma \ref{estbadset} implies that
\begin{equation*}
\begin{split}
\alpha:=\lim\limits_{\sigma\to 0+}\frac{\sigma E(\sigma)}{I(\sigma, Q_\sigma)}&=\lim\limits_{\sigma\to 0+}\Big(\frac{\sigma F(\sigma)}{I(\sigma, Q_\sigma)}\frac{E(\sigma)}{F(\sigma)}\Big)\\
&=\lim\limits_{\sigma\to 0+}\frac{\sigma F(\sigma)}{I(\sigma, Q_{\sigma})}\lim\limits_{\sigma\mapsto 0+}\frac{E(\sigma)}{F(\sigma)}<\infty.
\end{split}
\end{equation*}

\end{proof}

\begin{definition}
The value $\alpha$ is the \emph{order of $f$ at $0$} and denoted by $\alpha = \mathrm{ord}^f(0)$.
\end{definition}

\subsection{Improved monotonicity}

Using the H\"older regularity of Theorem \ref{chen} and the definition of $\alpha$ given by  Proposition \ref{orderdef}, we improve the H\"older result to have exponent corresponding to the order function. Such a result allows us to immediately conclude Lipschitz regularity whenever $\alpha \geq 1$.

\begin{lemma}\label{tildeEmonotone}Let  ${\bf B}$ be a local model, $g$  a normalized
Lipschitz metric defined on ${\bf B}(r)$, $(Y,d)$ a CAT(1) space and $f:({\bf B}(r),g) \rightarrow \B_\tau(P) \subset Y$  an energy minimizing map where $0<\tau< \frac \pi 4$. 
Let $\alpha= \mathrm{ord}^f(0)$ and $\gamma>0$ be the H\"older exponent of Theorem \ref{chen}. There exist constants $c_0, c_0'$ and $\sigma_0$ depending only on ${\bf B}(r)$, $E^f$, and the Lipschitz bound and the ellipticity constant of $g$ so that if $C$ is as in \eqref{e08} and
\[
\tilde E(\sigma) := E(\sigma)\exp\left(c_0\int_{A \cap(0,\sigma)} Cs^{2\gamma} \frac{E'(s)}{E(s)} ds \right),
\]then
\[
\sigma \mapsto e^{c_0 \sigma+c_0'\sigma^{2\gamma}} \frac{\tilde E(\sigma)}{\sigma^{2\alpha + n -2}} 
\]is non-decreasing for $\sigma \in (0, \sigma_0)$.
\end{lemma}
\begin{proof}

Let
\[
G(\sigma) = E(\sigma) \exp\left(\int_{A \cap(0,\sigma)} Cs^{2\gamma}\frac{E'(s)}{E(s)} ds \right).
\]Since \eqref{e010}, \eqref{e019} imply that $\exp \left(\int_{A \cap(0,1)} Cs^{2\gamma}\frac{E'(s)}{E(s)} ds \right)$ is finite and by definition
\[
G(\sigma) = F(\sigma) \exp \left(\int_{A \cap(0,1)} Cs^{2\gamma} \frac{E'(s)}{E(s)} ds \right), 
\]the monotonicity of \eqref{Fnondec} implies that
\[
\sigma \mapsto e^{c_3\sigma}\frac{\sigma G(\sigma)}{I(\sigma)}
\]is non-decreasing. Since 
$\int_{A \cap (0,\sigma)}Cs^{2\gamma}\frac {E'(s)}{E(s)}ds \geq 0$, $E(\sigma) \leq G(\sigma)$.  Therefore,
\begin{equation}
 \alpha = \lim_{\sigma \rightarrow 0} \frac{\sigma
E(\sigma)}{I(\sigma, Q_{\sigma})} \leq \lim_{\sigma \rightarrow 0} \frac{\sigma
G(\sigma)}{I(\sigma, Q_{\sigma})}=\lim_{\sigma \rightarrow 0}
e^{c_3\sigma} \frac{\sigma G(\sigma)}{I(\sigma, Q_{\sigma})}.
\end{equation} Thus for $\sigma$ sufficiently small,
\begin{equation} \label{uk}
\alpha \leq   e^{c_3 \sigma} \frac{\sigma G(\sigma)}{I(\sigma, Q_{\sigma})}= e^{c_3 \sigma} \frac{\sigma E(\sigma)}{I(\sigma, Q_{\sigma})} \cdot \exp\left(  \int_{A \cap (0,\sigma)} Cs^{2\gamma}\frac{E'(s)}{E(s)} ds \right)
\end{equation}
and 
\begin{equation} \label{uK}
\frac{\sigma E(\sigma)}{I(\sigma, Q_{\sigma})} \leq e^{c_3 \sigma} \frac{\sigma G(\sigma)}{I(\sigma, Q_{\sigma})} \leq e^{c_3} \frac{G(1)}{I(1)}=:K.
\end{equation}
We now use arguments from the proof of Lemma \ref{estbadset}. By \eqref{Jest}, for $n$ such that
$\theta_0^{n+1} \leq \sigma < \theta_0^{ n}$,
\begin{eqnarray}
\int_{A \cap (0,\sigma)} Cs^{2\gamma}\frac{E'(s)}{E(s)} ds & \leq & \int_0^{\sigma}Cs^{2\gamma} \frac{J'(s)}{J(s)} ds \nonumber \\
& \leq & \sum_{k=n}^{\infty} \int_{\theta_0^{k+1}}^{\theta_0^{k}} Cs^{2\gamma} \frac{J'(s)}{J(s)} ds \nonumber \\
& \leq & \left(\sum_{k=n}^{\infty}C\theta_0^{ 2\gamma k}\right)\log \theta_0^{-j } \nonumber \\
&=&\left(\frac{C\theta_0^{  2\gamma n}}{1-\theta_0^{2 \gamma}} \right) \log\theta_0^{-j} \nonumber \\
& \leq & c_4\theta_0^{2 \gamma n}\nonumber \\
& \leq & \frac{c_4}{\theta_0^{2\gamma}} \sigma^{2\gamma} =: c_5 \sigma^{2\gamma}. \label{linear}
\end{eqnarray}

By \eqref{uk} and \eqref{linear}
\begin{equation}\label{alphaineq}
\alpha \leq e^{c_3\sigma +c_5\sigma^{2\gamma}} \frac{\sigma E(\sigma)}{I(\sigma,Q_\sigma)}\leq e^{c_6\sigma^{2\gamma}}\frac{\sigma E(\sigma)}{I(\sigma,Q_\sigma)}. 
\end{equation}
By  \eqref{e03} and \eqref{EIinequ},
\begin{align*}
 2\Big(1-c'\sigma^{2\gamma}\Big)E(\sigma) &\leq  \int_{\partial {\bf B}(\sigma)} \frac \partial{\partial r} d^2(f,f(0)) d\Sigma_g + I(\sigma) +k\sigma^2E'(\sigma)\\
 & \leq I'(\sigma) -\frac{n-1}\sigma I(\sigma) + (1+c_2)I(\sigma) + k\sigma^2 E'(\sigma).
\end{align*}Here $k, c_2$ depend on the Lipschitz bound of $g$. Therefore
\begin{align*}
\frac{\frac{2\left(1-c'\sigma^{2\gamma}\right)\sigma E(\sigma)}{I(\sigma)} + n-1 -O(\sigma)}\sigma &\leq \frac{I'(\sigma)}{I(\sigma)} + k\sigma^2 \frac{E'(\sigma)}{I(\sigma)}\\
& \leq \frac{I'(\sigma)}{I(\sigma)} + kK \sigma\frac{E'(\sigma)}{E(\sigma)} \text{ by } \eqref{uK}\\
& \leq \frac{G'(\sigma)}{G(\sigma)}+\frac 1\sigma + kK \sigma\frac{E'(\sigma)}{E(\sigma)}+c_3.
\end{align*}For the last inequality we use \eqref{e08} and \eqref{diffF} to show that
\[
\frac{G'(\sigma)}{G(\sigma)} = \frac{F'(\sigma)}{F(\sigma)} \geq \frac{I'(\sigma)}{I(\sigma)} - \frac 1\sigma-c_3.
\]Note further that by applying \eqref{alphaineq} and absorbing the higher order terms of the exponential into $O(\sigma)$, there exists $c_0'$ such that
\[
\frac{G'(\sigma)}{G(\sigma)} + kK \sigma\frac{E'(\sigma)}{E(\sigma)}\geq \frac{2\alpha-c_0'\sigma^{2\gamma}}{\sigma} + \frac{n-2-O(\sigma)}{\sigma}   - c_3.
\]
From this point forward, we presume that $\sigma \in A$. Indeed, if $\sigma \notin A$ then the appropriate differential inequality is satisfied which immediately proves the monotonicity. By definition
\[
\frac{\tilde E'(\sigma)}{\tilde E(\sigma)} \geq \frac{E'(\sigma)}{E(\sigma)}\left(1+ c_0C\sigma^{2\gamma}\right).
\]Choose $c_0$ sufficiently large so that $c_0C \geq kK+C$. Then, since $\sigma \in A$, 
\[
\frac{\tilde E'(\sigma)}{\tilde E(\sigma)} \geq \frac{E'(\sigma)}{E(\sigma)}\left(1+ (kK+C)\sigma^{2\gamma}\right) \geq \frac{G'(\sigma)}{G(\sigma)}+ kK\sigma\frac{E'(\sigma)}{E(\sigma)}.
\]Then, we may increase $c_0$ if necessary to determine that
\begin{align*}
 \frac{\tilde E'(\sigma)}{\tilde E(\sigma)}  
 &\geq \frac{n-2+2\alpha}{\sigma} -c_0'\sigma^{2\gamma-1}-c_0.
\end{align*}

It follows that
\[
\frac d{d\sigma} \log \left( \frac{e^{c_0\sigma+ c_0'\sigma^{2\gamma}} \tilde E(\sigma)}{\sigma^{n-2+2\alpha}}\right) 
 \geq 0.\]
\end{proof}

\begin{corollary}\label{alphamonotonecor}
Let ${\bf B}$ be a dimension-$n$, codimension-$\nu$ local model, $g$ a normalized Lipschitz metric  defined on ${\bf B}(r)$, and $(Y,d)$ a CAT(1) space and $f:({\bf B}(r),g) \rightarrow \B_\tau(P) \subset Y$  an energy minimizing map where $0<\tau< \frac \pi 4$. Let $\alpha=\mathrm{ord}^f(0)$ and $\gamma>0$ be the H\"older exponent of Theorem \ref{chen}. Then there exist constants $c,k>0$, $\sigma_0<1$ depending on ${\bf B}(r)$, $E^f$, and the Lipschitz bound and the ellipticity constant of $g$ so that
\begin{equation}\label{monotoneorder}
\frac{\sigma E(\sigma)}{I(\sigma,Q_\sigma)} \leq ce^{k \rho}\frac{\rho E(\rho)}{I(\rho,Q_\rho)} \text{ for } 0<\sigma \leq\rho \leq \sigma_0,
\end{equation}and
\begin{equation}\label{Emonotoneorder}
\frac{E(\sigma)}{\sigma^{n-2+2\alpha}} \leq e^{k\rho^{2\gamma}}\frac{E(\rho)}{\rho^{n-2+2\alpha}} \text{ for } 0<\sigma \leq\rho \leq \sigma_0.
\end{equation}
\end{corollary}

\begin{proof}
For $F(\sigma)$ defined as in \eqref{Fdef}, \eqref{e010} and \eqref{e019} imply that
\[
\frac 1cE(\sigma) \leq F(\sigma) \leq E(\sigma)
\]for some $c>1$ depending on the H\"older constant from Theorem \ref{chen}. Using this uniform bound and the fact that $ e^{c_3\sigma}\frac{\sigma F(\sigma)}{I(\sigma, Q_\sigma)}$ is monotone by \eqref{Fnondec} implies \eqref{monotoneorder}.

Lemma \ref{tildeEmonotone}, the definition of $\tilde E$, and the fact that \eqref{linear} gives the bound
\begin{equation*}\label{EtildeE}
E(\sigma) \leq \tilde E(\sigma) \leq e^{c_0c_5\sigma^{2\gamma}}E(\sigma).
\end{equation*} together imply \eqref{Emonotoneorder}.
\end{proof}

As in the conclusion of Section \ref{monotoness}, we now consider monotonicity for metrics $g$ that are not necessarily normalized. Recall that if $g$ is a Lipschitz metric and $h:= L^*_xg$, where $L_x$ is the map given by \cite[Proposition 2.1]{daskal-meseCAG}, then $h$ is normalized. Moreover, when $f$ is minimizing with respect to the metric $g$, then $f \circ L_x$ is minimizing with respect to $h$.
\begin{definition}\label{orderdef2}
For $f$ minimizing with respect to a Lipschitz metric $g$, we define the \emph{order of $f$ at $x$} as
\[
\alpha_x= \mathrm{ord}^f(x) := \mathrm{ord}^{f \circ L_x}(0).
\]
\end{definition}Recalling that
\[
E_x(\sigma):= \int_{B_x(\sigma)}|\nabla f|^2 d\mu_g,
\]and $B_x(\sigma)$ denotes the Euclidean ball about $x$ of radius $\sigma$, 
we prove the monotonicity of Proposition \ref{relax} with exponent $n-2+2\alpha$.
\begin{proposition} \label{relax2}
Let ${\bf B}$ be a dimension-$n$, codimension-$\nu$ local model, $g$ a  Lipschitz metric  defined on ${\bf B}(r)$ with ellipticity constant $\lambda \in (0,1]$, $(Y,d)$ a CAT(1) space and $f:({\bf B}(r),g) \rightarrow \B_\tau(P) \subset Y$  an energy minimizing map where $0<\tau< \frac \pi 4$. Then there exists $C \geq 1$ depending on ${\bf B}(r)$, the  Lipschitz bound and the ellipticity constant of $g$ so that for every $x \in {\bf B}(r)$, 
\begin{equation} \label{gooddecay2}
\frac{E_x(\sigma)}{ \sigma^{n-2+2\alpha_x}} \leq C \frac{E_x(\rho)}{ \rho^{n-2+2\alpha_x}} , \ 0<\sigma<\rho\leq   r(x)
\end{equation}
where $r(x)$ is defined as in \eqref{defr(x)}.
\end{proposition}
\begin{proof}
Following the proof of Proposition \ref{relax}, the result follows from \cite[Proposition 2.1]{daskal-meseCAG} and Corollary \ref{alphamonotonecor}. 
\end{proof}

Using \eqref{gooddecay2} and the techniques of Section \ref{Holderss}, we immediately determine H\"older regularity for $f$ with exponent depending on the order function.

\begin{theorem} \label{chen2}
Let  ${\bf B}$ be a local model, $g$  a 
Lipschitz Riemannian metric defined on ${\bf B}(r)$, $(Y,d)$ a CAT(1) space and $f:({\bf B}(r),g) \rightarrow \B_\tau(P) \subset Y$  an energy minimizing map where $0<\tau< \frac \pi 4$.  If $0<\alpha \leq \alpha_x$ for all $ x\in {\bf B}(\varrho r)$ where $\varrho \in(0,1)$,  then  there exists $C$ depending only on the Lipchitz bound and ellipticity constant of $g$, $E^f$, ${\bf B}(r)$ and $\varrho$ such  that
\[
d(f(x),f(y)) \leq C |x-y|^{\alpha}, \hspace{0.1in} \forall x,y \in
{\bf B}(\varrho r).
\]
\end{theorem}

\begin{proof}
The result follows immediately from \eqref{gooddecay2},  Proposition~\ref{propjoe}, the Poincar\'e inequality of \cite[Theorem 2.7]{daskal-meseCAG}, and Lemma~\ref{campanato}.
\end{proof}

\section{Tangent Map Construction}\label{Tangentss}
Given a domain $\Omega$, NPC spaces $(Y_k,d_k)$,  and maps $f_k:\Omega \to Y_k$, Korevaar and Schoen \cite[Section 3]{korevaar-schoen2} develop the notion of convergence of maps in the pullback sense. This allowed \cite{daskal-meseCAG} to define a tangent map of $f:{\bf B}(r) \to Y$ when $Y$ is NPC. They then related the homogeneity of a tangent map to the order of $f$ and used this to get the Lipschitz regularity.

Rather than reconstruct the entire argument when $Y$ is CAT(1), we will consider the {tangent map} of $f$ that is determined by the tangent map construction in \cite{daskal-meseCAG} for the lifted map $\overline f:{\bf B}(r) \to \mathcal CY$. Since $\mathcal CY$ is NPC, we do not need to reconstruct the theory. Instead, we use the minimizing property of $f$ to prove that the proposed tangent map exists. 

\subsection{Limit maps in the pullback sense}
We first recall the construction in \cite[Section 3]{korevaar-schoen2} and its extension to local models in \cite[Section 5]{daskal-meseCAG}.

Let $f_k: {\bf B}(r) \to Y_k$ where each $(Y_k,d_k)$ is an NPC space. Since each $f_k$ maps to a different metric space, convergence cannot be understood in a pointwise sense without further work. If one considers the closed convex hull of each set $f_k({\bf B}(r))$ and corresponding pseudodistances $d_{k,\infty}$, convergence can be well understood by considering convergence of the pseudodistances. The construction proceeds as follows. 

Let $f:{\bf B}(r) \to Y$ and denote $\Omega_0={\bf B}(r)$, $f_0=f$, and let $d_0:\Omega_0 \times \Omega_0 \to \mathbb R^+ \cup \{0\}$ be the pseudodistance function $d_0(x,y):= d(f_0(x), f_0(y))$. Inductively define $\Omega_{i+1} = \Omega_i \times \Omega_i \times [0,1]$ and identify $\Omega_i \subset \Omega_{i+1}$ via the inclusion $x \mapsto (x,x,0)$. Define $f_{i+1}:\Omega_{i+1} \to Y$ by 
\[
f_{i+1}(x,y,t) = (1-t)f_i(x) + t f_i(y).
\]Let
\[
d_{i+1}(x,y):= d(f_{i+1}(x),f_{i+1}(y)).
\]Then
\[
d_{i+1}((x,x,0),(y,y,0)) = d_i(x,y),
\]
\[
d_{i+1}((x,y,s), (x,y,t)) = |s - t|d_i(x,y),
\]
\begin{align*}
d_{i+1}^2(z,(x,y,s))\leq &(1-s)d^2_{i+1}(z,(x,x,0)) + s d^2_{i+1}(z,(y,y,0))\\
& \quad\quad -s(1-s)d^2_{i+1}((x,x,0),(y,y,0)).
\end{align*}

Set $\Omega_\infty = \cup_i \Omega_i$ and define $f_\infty:\Omega_\infty \to (Y,d)$ such that $f_\infty|_{\Omega_i} = f_i$. Define a pseudodistance function $d_\infty(x,y):= d(f_\infty(x), f_\infty(y))$. Define the metric space $(Y_*, d_*)$ as the completion of the quotient metric space constructed from $(\Omega_\infty, d_\infty)$, where the quotient space is defined via the equivalence relation of zero pseudodistance.
Then by construction, in particular the properties of $d_{i}$, $(Y_*,d_*)$ is an NPC space. 
Moreover, $(Y_*,d_*)$ is isometric to $\mathrm{Cvx}(f({\bf B}(r)))$, the closed convex hull of $f({\bf B}(r))$, with isometry given by the unique extension of $f_\infty$ to $Y_*$.

\begin{definition}  \label{pullbacksense}
Let $v_k:{\bf B}(r) \rightarrow (Y_k,d_k)$ be a sequence of maps
to NPC spaces.  We say $v_k$ \emph{converges to $v_*$ in the pullback
sense} if the corresponding pullback pseudodistances $d_{k,\infty}$ converge pointwise to a limit pseudodistance $d_{\infty}$ on $\Omega_\infty \times \Omega_\infty$, and $v_*= \pi \circ \iota$ where $\iota:{\bf B}(r) \to \Omega_\infty$ is the inclusion map and $\pi$ is the natural projection map of $\Omega_\infty$ onto the metric completion $(Y_*,d_*)$ of the quotient space constructed from $(\Omega_\infty, d_\infty)$. 
\end{definition}

Given $v_*$ as above, we can replace $f$ in the outlined construction by $v_*$. Then $d_{*,i}$ denotes the corresponding pullback pseudodistance function of $v_{*,i}$ and $d_{*,\infty}$ denotes the corresponding pullback pseudodistance function of $v_{*,\infty}$. 
In this case, $(d_*)_*=d_{*}$.

\begin{definition}  Suppose $v_k$ converge to $v_*$ in the pullback
sense.  Let $d_{k,i}$ (resp. $d_{\infty,i}$) be the corresponding
pullback pseudodistance function to $v_{k,i}:\Omega_i \rightarrow
(Y_k,d_k)$ (resp. $v_{*,i}: \Omega_i \rightarrow (Y_*,d_*)$).  We
say that the convergence is \emph{locally uniform} if the convergence of
$d_{k,i}$ to the limit $d_{*,i}$  is uniform on each compact
subset of $\Omega_i \times \Omega_i$.  In this case, we also say
 $v_k \to v_*$ \emph{locally uniformly in the pullback sense}.
\end{definition}

\begin{proposition}\label{MOC}\cite[Proposition 3.7]{korevaar-schoen2}, \cite[Proposition 5.1]{daskal-meseCAG}
Let $v_k:{\bf B}(r) \to (Y_k,d_k)$ be a sequence of maps to NPC spaces $Y_k$ for which there is uniform modulus of continuity control. That is, assume that for each $x \in {\bf B}(r)$ and $R>0$ there exists a positive function $\omega(x,R)$ which is monotone in $R$ satisfying
\[
\lim_{R \to 0} \omega(x,R) = 0,
\] and so that for each $k \in \mathbb Z$
\[
\max_{y \in B(x,R)} d_k(v_k(x), v_k(y)) \leq \omega(x,R).
\]Then there is a subsequence of $v_k$ which converges locally uniformly in the pullback sense to a limit map $v_*:{\bf B}(r) \to (Y_*,d_*)$. Moreover, $v_*$ satisfies the same modulus of continuity estimates.
\end{proposition}
\subsection{The tangent map construction}
Let ${\bf B}$ be a dimension-$n$, codimension-$\nu$ local model and $g$ a normalized Lipschitz metric on ${\bf B}(1)$. For $r \in (0,1)$ and a map $f:{\bf B}(r) \to Y$ we will consider the \emph{$\lambda$-blow up map} $f_\lambda:{\bf B}(r/\lambda) \to (Y, d_\lambda)$ and the \emph{lifted $\lambda$-blow up map} $\overline f_\lambda:{\bf B}(r/\lambda) \to (\mathcal CY, D_\lambda)$ where
\begin{align*}
g_\lambda(x)&:= g(\lambda x)\\
\mu_\lambda^d&:= (\lambda^{1-n}\, \,^dI(\lambda))^{1/2}\\
\mu_\lambda^D&:= (\lambda^{1-n}\, \,^DI(\lambda))^{1/2}\\
d_\lambda(P,Q)&:= (\mu_\lambda^d)^{-1}d(P,Q)\\
D_\lambda(P,Q)&:=(\mu_\lambda^D)^{-1}D(P,Q)\\
f_\lambda(x)&:=f(\lambda x)\in Y\\
\overline f_\lambda(x)&:= [f(\lambda x),1] \in Y \times \{1\} \subset \mathcal CY.
\end{align*}Above we have denoted
\begin{align*}^dI(\lambda)&:= \inf_{Q \in Y} \int_{\partial {\bf B}(\lambda)} d^2(f,Q) d\Sigma_g \\^DI(\lambda)&:= \inf_{Q \in \mathcal CY} \int_{\partial {\bf B}(\lambda)} D^2([f,1],Q) d\Sigma_g.
\end{align*}

\begin{definition}
If there exists an NPC space $(Y_*,d_*)$ and a sequence $\lambda_k \to 0$ such that $f_{\lambda_k}$ converges locally uniformly in the pullback sense to $\overline f_*:{\bf B} \to (Y_*,d_*)$ then $\overline f_*$ is called a \emph{tangent map of $f$}.
\end{definition}

\begin{proposition}\label{tangentmapprop}
For ${\bf B}$ a dimension-$n$, codimension-$\nu$ local model, $g$ a normalized Lipschitz metric defined on ${\bf B}(r)$, and $(Y,d)$ a CAT(1) space, let $f:({\bf B}(r),g) \to \B_\tau(P) \subset Y$ be an energy minimizing map where $0<\tau< \frac \pi 4$. Then $f$ has a tangent map $\overline f_*:({\bf B},g) \to (Y_*,d_*)$. 
Moreover, $\overline f_*:{\bf B}(1) \to (Y_*,d_*)$ is a non-constant, energy minimizing map. 
\end{proposition}
\begin{proof}
We first determine uniform modulus of continuity control on the maps $\overline f_\lambda$. Following \cite[Lemma 6.1]{daskal-meseCAG}, for the maps $f_\lambda:{\bf B}(r) \to (Y, d_\lambda)$, for sufficiently small $\lambda>0$
\[
^{d_\lambda}E^{f_\lambda}_{g_\lambda}[{\bf B}(1)] = \int_{{\bf B}(1)} |\nabla f_\lambda|^2 d\mu_{g_\lambda} \leq 2 \alpha
\]where $\alpha$ is the order of $f$ at zero. Given the uniform Lipschitz bounds on the metrics $g_\lambda$, we appeal to Theorem \ref{chen} and note that
\begin{equation}\label{unifMOC}
d_\lambda(f_\lambda(x),f_\lambda(y)) \leq C|x-y|^\gamma \text{ for all } x,y \in {\bf B}(r),
\end{equation}where $C,\gamma$ are independent of $\lambda$. This immediately implies uniform modulus of continuity control on the maps $f_\lambda$ but not on their lifted maps $\overline f_\lambda$. To determine the necessary control for the lifted maps, we will consider the relation between $d_\lambda$ and $D_\lambda$ as $\lambda \to 0$.

Since $f$ is energy minimizing into $Y$, Theorem \ref{chen} implies that $f(\partial {\bf B}(\lambda)) \subset \B_{C\lambda^\gamma}(f(0)) \subset Y$ where $C$ depends only on the Lipschitz bound and ellipticity constant of $g$, $E^f$, and ${\bf B}(r)$. By \eqref{D_d_compare}, given $Q \in \B_{C\lambda^\gamma}(f(0))$, for all $x \in \partial {\bf B}(\lambda)$,
\begin{equation}\label{lambdadistcompare}
D^2(\overline f(x),[Q,1]) \leq d^2(f(x),Q).
\end{equation}It follows that 
for $Q_\lambda^d$ such that $^dI(\lambda,Q_\lambda^d) = \inf \,\,^dI(\lambda, Q)$,
\[
^DI(\lambda) \leq \,\,^DI(\lambda,[Q_\lambda^d,1]) \leq \,\,^dI(\lambda).
\]Let $Q_\lambda^D \in \mathcal CY$ such that $^DI(\lambda,Q_\lambda^D) = \inf \,\,^DI(\lambda, Q)$. Then $Q_\lambda^D \in \B_{C\lambda^\gamma}([f(0),1]) \subset \mathcal CY$ and by
\eqref{D_d_compare2},
\[
^dI(\lambda) \leq \,\,^dI(\lambda,\pi_1(Q^D_\lambda)) \leq (1+C\lambda^{2\gamma})\,\,^DI(\lambda),
\]where $\pi_1:\mathcal CY \to Y$ is the projection map onto the first component of $Y \times [0,\infty)$. Therefore, for $\mu_\lambda^d, \mu_\lambda^D$ the rescalings of $d,D$ respectively, and $\lambda>0$ sufficiently small,
\begin{equation}\label{mucomparison}
\frac 12 \leq (1+C\lambda^{2\gamma})^{-1/2} \leq \frac{\mu_\lambda^D}{\mu_\lambda^d} \leq 1.
\end{equation}It follows by \eqref{D_d_compare}, \eqref{lambdadistcompare}, and \eqref{mucomparison} that for sufficiently small $\lambda>0$ and all $x,y \in {\bf B}(r)$,
\begin{align*}
D_\lambda(\overline f_\lambda(x),\overline f_\lambda(y))&= (\mu_\lambda^D)^{-1}D([f(\lambda x),1], [f(\lambda y),1])\\
&\leq 2(\mu_\lambda^d)^{-1}D([f(\lambda x),1], [f(\lambda y),1])\\
&\leq 2(\mu_\lambda^d)^{-1}d(f(\lambda x), f(\lambda y)) \\
&=2 d_\lambda(f_\lambda(x),f_\lambda(y)).
\end{align*}By \eqref{unifMOC}, the maps $\overline f_\lambda$ into $(\mathcal CY, D_\lambda)$ possess uniform modulus of continuity control. Therefore, by Proposition \ref{MOC}, there exists a sequence $\lambda_k \to 0$ and an NPC space $((\mathcal CY)_*,D_*)$ such that $\overline f_{\lambda_k}$ converge locally uniformly in the pullback sense to a limit map $\overline f_*:{\bf B}(1) \to ((\mathcal CY)_*,D_*)$.
\begin{claim}
$\overline f_*$ is a tangent map of $f$.
\end{claim} 
\begin{proof}
We need to show that $d_{\lambda,n} \to d_{*,n}$ uniformly on 
$\Omega_n\times \Omega_n$ for all $n \in \mathbb N \cup \{0\}$. Since the uniform convergence for $D_{\lambda,n} \to d_{*,n}$ is already 
established, and since $\mu^d_\lambda/\mu^D_\lambda \to 1$ uniformly by \eqref{mucomparison}, it is enough to show that for all $n \in \mathbb N \cup \{0\}$, 
 \begin{equation}\label{dist_i_conv}
\frac {d(f_{\lambda,n}(x), f_{\lambda,n}(y))}{D(\overline f_{\lambda,n}(x),\overline f_{\lambda,n}(y))} \to 1
 \end{equation} uniformly for all $x\neq y \in \Omega_n$. Proceeding by induction requires that we also demonstrate that
 \begin{equation}\label{DistProjUnif}
D([f_{\lambda,n}( x),1], \overline f_{\lambda, n}( x)) \to 0
\end{equation} uniformly for $ x \in \Omega_{n}$.

Observe that by \eqref{projectiondistance},
\[
\frac{d(f_\lambda(x),f_\lambda(y))}{D(\overline f_\lambda(x), \overline f_\lambda(y))}  \to 1 
\]uniformly for all $x\neq y \in {\bf B}(1)= \Omega_0$ so \eqref{dist_i_conv} holds easily for $n=0$. Morover, \eqref{DistProjUnif} is trivial for $n=0$ since $\overline f_{\lambda,0}(x)= [f_{\lambda,0}(x),1] \in \mathcal CY$.

Now suppose that 
\begin{equation*}\label{iminus1dist}
\frac{d(f_{\lambda, i-1}( x),f_{\lambda, i-1}( y))}{D(\overline f_{\lambda, i-1}( x),\overline f_{\lambda, i-1}(  y))} \to 1
\end{equation*}uniformly for $ x \neq y \in \Omega_{i-1}$ and that 
\begin{equation}\label{iminus1proj}
D([f_{\lambda,i-1}( x),1], \overline f_{\lambda, i-1}( x)) \to 0
\end{equation} uniformly for $ x \in \Omega_{i-1}$.   We claim that together these imply that \eqref{dist_i_conv} and \eqref{DistProjUnif} hold for $n=i$ and $\mathbf x \neq \mathbf y \in \Omega_{i}$. 

Consider $\mathbf x , \mathbf y \in \Omega_i$ with $\mathbf x = (x_1, x_2,s)$ and $\mathbf x \neq \mathbf y$. Since, by Theorem \ref{chen}, $f_{\lambda,i}({\bf B}(1)) \subset \B_{C\lambda^\gamma}(f(0))$, \eqref{projectiondistance} implies that
 \[
 \frac{d(f_{\lambda,i}(\mathbf x), f_{\lambda,i}(\mathbf y))}{D([f_{\lambda,i}(\mathbf x),1], [f_{\lambda,i}(\mathbf y),1])} \to 1\text{ uniformly for } \mathbf x \neq \mathbf y \in \Omega_i.
 \] Thus, it is enough to show that 
 \[
 \frac{D([f_{\lambda,i}(\mathbf x),1], [f_{\lambda,i}(\mathbf y),1])}{D( \overline f_{\lambda,i}(\mathbf x),\overline f_{\lambda,i}(\mathbf y))}  \to 1 \text{ uniformly for } \mathbf x \neq \mathbf y \in \Omega_i.
 \]
 Note that if $x_1=x_2$ then $f_{\lambda,i}(\mathbf x) = f_{\lambda,i-1}(x_1)$ and $\overline f_{\lambda,i}(\mathbf x) = \overline f_{\lambda,i-1}(x_1)$. Thus $D([f_{\lambda,i}(\mathbf x),1], \overline f_{\lambda,i}(\mathbf x)) \to 0$ uniformly by \eqref{iminus1proj}. 
 Now suppose that $x_1 \neq x_2$. By hypothesis, with $\gamma_\lambda \subset Y$ the geodesic connecting $f_{\lambda,i-1}(x_1)$ to $f_{\lambda,i-1}(x_2)$ and $\overline \gamma_\lambda \subset \mathcal CY$ the geodesic connecting $\overline f_{\lambda,i-1}(x_1)$ to $\overline f_{\lambda,i-1}(x_2)$,
\begin{equation}\label{DistRatioUnif}
\frac{\ell(\gamma_\lambda)}{\ell(\overline \gamma_\lambda)} = \frac{d(f_{\lambda, i-1}(x_1),f_{\lambda, i-1}(x_2))}{D(\overline f_{\lambda, i-1}(x_1),\overline f_{\lambda, i-1}(x_2))} \to 1
\end{equation}uniformly. For $j=1,2$, let $\overline \gamma^j_\lambda \subset \mathcal CY$ be the geodesic connecting $\overline f_{\lambda,i-1}(x_j)$ to $[f_{\lambda, i}(\mathbf x),1]$. Then by the triangle inequality,
\[
\ell(\overline \gamma_\lambda) \leq \ell(\overline \gamma^1_\lambda) + \ell(\overline \gamma_\lambda^2) \leq \ell (\gamma_\lambda) + \sum_{j=1,2}D(\overline f_{\lambda, i-1}(x_j), [f_{\lambda,i-1}(x_j),1]).
\]Thus by \eqref{iminus1proj} and \eqref{DistRatioUnif},
\[
\frac{\ell(\overline \gamma^1_\lambda) + \ell(\overline \gamma_\lambda^2) }{ \ell (\overline \gamma_\lambda) } \to 1 \text{ uniformly}.
\]We consider the geodesic triangle in $\mathcal CY$ with endpoints $\overline f_{\lambda,i-1}(x_1), \overline f_{\lambda,i-1}(x_2), [f_{\lambda,i}(\mathbf x),1]$. Using a comparison triangle in $\mathbb R^2$, the side length relation implies that 
\[
D([f_{\lambda,i}(\mathbf x),1], \overline f_{\lambda,i}(\mathbf x)) \to 0\text{ uniformly for all } \mathbf x \in \Omega_i.
\] Therefore \eqref{DistProjUnif} holds for $n=i$. By the triangle inequality,
\begin{equation}\label{D_i_Unif}
D([f_{\lambda,1}(\mathbf x),1], [f_{\lambda,1}(\mathbf y),1]) - D(\overline f_{\lambda,1}(\mathbf x),\overline f_{\lambda,1}(\mathbf y))\to 0 \text{ uniformly}
\end{equation}and thus \eqref{dist_i_conv} holds for $n=i$ and $\mathbf x \neq \mathbf y$.

Therefore $f_{\lambda_k}$ converges uniformly locally in a pullback sense to $\overline f_*$ and it is reasonable to consider the target using the notation $(Y_*,d_*)$.

\end{proof}

Finally, we prove that $\overline f_*$ is minimizing. Let $v_{\lambda} = ^{Dir}\overline{f}_{\lambda}: {\bf B}(1) \to (\mathcal CY,D_\lambda)$ denote the Dirichlet solution for $\overline f_\lambda$. 
As before, we note that $f_{\lambda}( {\bf B}(1))\subset \B_{C\lambda^{\gamma}}(f(0))$, and so 
by \eqref{D_d_compare}, it follows that $\overline{f}_{\lambda}( {\bf B}(1)) \subset \B^{\mathcal{C}Y}_{C\lambda^{\gamma}}([f(0),1])$.
Now, since $\mathcal CY$ is an NPC space, projection onto convex domains decreases energy. Therefore, since $v_{\lambda}$ is an energy minimizer and $v_{\lambda}|_{\partial  {\bf B}(1)} = \overline{f}_{\lambda}|_{\partial  {\bf B}(1)}\subset \B^{\mathcal{C}Y}_{C\lambda^{\gamma}}([f(0),1])$,   it follows that $v_{\lambda}({\bf B}(1))\subset \B^{\mathcal{C}Y}_{C\lambda^{\gamma}}([f(0),1])$. Thus, $\pi_2(v_{\lambda})\geq 1-C\lambda^{\gamma}$, where $\pi_2:\mathcal{C}Y\rightarrow [0,\infty)$ is the projection onto the second component of the cone over $Y$.
For $\lambda$ sufficiently small, we apply \eqref{distancelowerbnd} 
 and observe that
\begin{equation*}
D^2(v_{\lambda}(x), v_{\lambda}(y))\geq (1-C\lambda^{\gamma})^2D^2(\Pi(v_{\lambda}(x)), \Pi(v_{\lambda}(y))),
\end{equation*}
and thus $^DE^{\Pi(v_{\lambda})}[ {\bf B}(1)]\leq (1-C\lambda^{\gamma})^{-2}\,{}^DE^{v_{\lambda}}[ {\bf B}(1)]$. By \eqref{liftenergy}, if $\pi_1:\mathcal CY \to Y$ denotes the projection onto the first component of the cone over $Y$, then 
$^dE^{\pi_1( v_{\lambda})}[ {\bf B}(1)] = {}^DE^{\Pi(v_{\lambda})}[ {\bf B}(1)]$, and thus 
\begin{equation}\label{DE1}
{}^dE^{\pi_1( v_{\lambda})}[ {\bf B}(1)]\leq (1-C\lambda^{\gamma})^{-2}\,{}^DE^{v_{\lambda}}[ {\bf B}(1)].
\end{equation}

Again using \eqref{liftenergy}, and noting that $f_{\lambda}$ is energy minimizing with respect to $d_{\lambda}$,  
\begin{align}\label{DE2}
{}^D E^{\overline{f}_{\lambda}}[ {\bf B}(1)]={}^{d}E^{f_{\lambda}}[ {\bf B}(1)]  &= (\mu_{\lambda}^d)^2\,{}^{d_{\lambda}}E^{f_{\lambda}}[ {\bf B}(1)] \leq (\mu_{\lambda}^d)^2\,{}^{d_{\lambda}}E^{\pi_1( v_{\lambda})}[ {\bf B}(1)] ={}^{d}E^{\pi_1( v_{\lambda})}[ {\bf B}(1)].
\end{align}
Combining \eqref{DE1} and \eqref{DE2} we observe that 
\[
 {}^D E^{\overline{f}_{\lambda}}[ {\bf B}(1)]\leq (1-C\lambda^{\gamma})^{-2}\,{}^DE^{v_{\lambda}}[ {\bf B}(1)]
 \] and therefore  
\begin{align*}
 {}^{D_{\lambda}} E^{\overline{f}_{\lambda}}[ {\bf B}(1)] = (\mu_{\lambda}^D)^{-2}\,{}^D E^{\overline{f}_{\lambda}}[ {\bf B}(1)]&\leq (\mu_{\lambda}^D)^{-2}(1-C\lambda^{\gamma})^{-2}\,{}^DE^{v_{\lambda}}[ {\bf B}(1)]\\&= (1-C\lambda^{\gamma})^{-2}\,{}^{D_{\lambda}}E^{v_{\lambda}}[ {\bf B}(1)].
 \end{align*}
  Finally, since $v_{\lambda}$ is energy minimizing with respect to $D_{\lambda}$, we have that
\begin{equation*}
{}^{D_{\lambda}}E^{v_{\lambda}}[ {\bf B}(1)]\leq {}^{D_{\lambda}} E^{\overline{f}_{\lambda}}[ {\bf B}(1)] \leq (1-C\lambda^{\gamma})^{-2}\,{}^{D_{\lambda}}E^{v_{\lambda}}[ {\bf B}(1)],
\end{equation*}
and so it follows from \cite[Theorem 3.11]{korevaar-schoen2} 
 that $\overline{f}_*$ is minimizing. The non-constancy of $\overline f_*$ follows exactly as in the proof of \cite[Proposition 3.3]{gromov-schoen}.

\end{proof}

\section{Higher Regularity Results}\label{Lipss}

\subsection{Lipschitz regularity}
The Lipschitz regularity of $f$ at points in $X-X^{(n-2)}$ will follow from regularity results for minimizing maps into an NPC space, once we show that the order of $f$ is bounded below by the order of its tangent map $\overline f_*$.
\begin{lemma}\label{alphaeq}
Let ${\bf B}$ be a dimension-$n$, codimension-$\nu$ local model, $g$ a normalized Lipschitz metric defined on ${\bf B}(r)$ and $(Y,d)$ a CAT(1) space. If $f:({\bf B}(r),g) \to \B_\tau(P) \subset Y$ is an energy minimizing map with $0<\tau < \frac \pi 4$, let $\overline f_*:{\bf B} \to (Y_*,d_*)$ denote a tangent map of $f$, constructed as in Proposition \ref{tangentmapprop}. Then
\[
\mathrm{ord}^{\overline f_*}(0) \leq \mathrm{ord}^f(0),  \quad  \text{ i.e., } \alpha_* \leq \alpha.
\]
\end{lemma}
\begin{proof}
Note that $\overline f_*$ is a minimizing map into an NPC space and thus by \cite[Corollary 3.1]{daskal-meseCAG}, \[
\alpha_* := \lim_{\sigma \to 0} \frac{\sigma \,{^{d_*}}E^{\overline f_*}(\sigma)}{{^{d_*}} I^{\overline f_*}(\sigma)} <\infty.
\]
Let $\lambda_k$ denote the sequence defining the lifted tangent map $\overline f$. Define $f_k:= f_{\lambda_k}$, $\mu_k:= \mu_{\lambda_k}^d$, $d_k := d_{\lambda_k}$, and $g_k:= g_{\lambda_k}$. Then,
\begin{align*}
\lim_{k \to \infty} \frac{\sigma \,\,^{d_k}E^{f_k}_{g_k}(\sigma)}{\,\,^{d_k}I^{f_k}_{g_k}(\sigma)}&=
\lim_{k \to\infty} \frac{\sigma  \mu_k ^{-2} \lambda_k^{2-n}\, {^{d}}E^{f}_g(\lambda_k\sigma)}{\mu_k^{-2} \lambda_k^{1-n}\, {^{d}}I^{f}_{g}(\lambda_k\sigma)}\\
&=\lim_{k \to \infty} \frac{\lambda_k\sigma   \,{^{d}}E^{f}_g(\lambda_k\sigma)}{{^{d}}I^{f}_{g}(\lambda_k \sigma)}\\
&=\alpha.
\end{align*}Therefore, it is enough to show that
\begin{align*}
\frac{\sigma  \,{^{d_*}}E^{\overline f_*}(\sigma)}{{^{d_*}}I^{\overline f_*}(\sigma)}\leq \lim_{k \to \infty} \frac{\sigma \,\,^{d_k}E^{f_k}_{g_k}(\sigma)}{\,\,^{d_k}I^{f_k}_{g_k}(\sigma)}.
\end{align*} 

By the pointwise convergence of  
$(\mu_{\lambda_k}^D)^{-1}D:=D_k \to D_*$ locally, uniformly on compact sets, 
\[^{D_k}I^{\overline f_k}_{g_k}(\sigma) \to \, {^{d_*}}I^{\overline f_*}(\sigma).
\]

By a change of variables and properties of the $\lambda_k$-blow up maps, 
\begin{align*}
^{d_k}I^{f_k}_{g_k}(\sigma) &= (\mu_k^d)^{-2} \lambda_k^{1-n}(^dI)^f_g(\lambda_k \sigma),\\
^{D_k}I^{\overline f_k}_{g_k}(\sigma) &= (\mu_k^D)^{-2} \lambda_k^{1-n}(^DI)^{\overline f}_g(\lambda_k \sigma).
\end{align*}The proof of Proposition \ref{tangentmapprop} immediately implies that 
\[
^{d_k}I^{f_k}_{g_k}(\sigma)-\,\,^{D_k}I^{\overline f_k}_{g_k}(\sigma) \to 0
\]as $k \to \infty$ and thus
\[
^{d_k}I^{f_k}_{g_k}(\sigma) \to  \, {^{d_*}}I^{\overline f_*}(\sigma) \text{ as } k \to \infty.
\]
By \cite[Lemma 8.1, Lemma 8.7]{korevaar-schoen2}, $\liminf_{k \to \infty} {^ {D_k}}E^{\overline f_k}_{g_k}(\sigma) \geq \, {^{d_*}}E^{\overline f_*}(\sigma)$. By \eqref{projectiondistance}, $^dE^{f_k}_{g_k}(\sigma) = \,^DE^{\overline f_k}_{g_k}(\sigma)$ for all $k$. Moreover, since $\mu_k^d - \mu_k^D \to 0$ as $k \to \infty$, 
$${^{d_*}}E^{\overline f_*}(\sigma) \leq \liminf_{k \to \infty} {^{d_k}}E^{f_k}_{g_k}(\sigma) .$$ This implies the result.
\end{proof}

\begin{proposition}
Let $\bf B$ be a dimension-$n$, codimension-$\nu$ local model with $\nu \in \{0,1\}$, $g$ a normalized Lipschitz metric defined on ${\bf B}(r)$, $(Y,d)$ a CAT(1) space and $f:({\bf B}(r),g) \to \B_\tau(P) \subset Y$ an energy minimizing map  with $0<\tau< \frac \pi 4$. Then $f$ is Lipschitz continuous in ${\bf B}(\varrho r)$ with Lipschitz constant depending on $\varrho \in (0,1)$, $({\bf B}(r),g)$, and the total energy of the map $f$.
\end{proposition}

\begin{proof}For each $x \in {\bf B}(\varrho r)$, consider the normalized map $f_x:= f \circ L_x$, minimizing with respect to the normalized metric $h=L^*_xg$. Here the map $L_x$ is as in \cite[Proposition 2.1]{daskal-meseCAG} but it plays a slightly different role. The metric on ${\bf B}$ is  normalized at the origin of ${\bf B}$. To determine a lower bound for the order of $f$ at $x$, we want to consider its tangent map about $x$, which requires that the metric be normalized about $x$. $L_x$ does this for us, and thus we may consider the tangent map of $f_x$ at $0$ (or $f$ at $x$).

By Proposition \ref{tangentmapprop}, the tangent map $(f_x)_*:{\bf B}(1) \to Y_*$ is minimizing into the NPC space $Y_*$. Therefore, by \cite[Lemma 8.1, Lemma 8.7]{daskal-meseCAG}, $\alpha_* \geq 1$. Lemma \ref{alphaeq} implies that $\mathrm{ord}^{f_x}(0) \geq 1$ and thus $\mathrm{ord}^f(x)\geq 1$. The result now follows from Theorem \ref{chen2}.
\end{proof}

The Lipschitz regularity, Theorem \ref{LipThm} item (1), follows immediately from \cite[Proposition 2.1]{daskal-meseCAG}.

\subsection{Regularity at a higher codimension singular point}
Given a Riemannian complex $X$ and an NPC space $T$, we define the center of mass of a map $u \in L^2(X,T)$ to be the unique point $\overline u \in T$ (with existence and uniqueness given by \cite[Proposition 2.5.4]{korevaar-schoen1}) such that
\[
\int_X d^2_T(u, \overline u) \, d\mu_g := \inf_{P \in T} \int_X d_T^2(u, P) \, d\mu_g.
\]We define the first eigenvalue of $X$ with values in $T$ by the Rayleigh quotient 
\[
\lambda_1(X,T):= \inf_{u \in W^{1,2}(X,T)} \frac{\int_X |\nabla u|^2 \, d\mu_g}{\int_X d^2(u, \overline u) \, d\mu_g}.
\]In application, $T$ will be the tangent cone of the CAT(1) space $Y$ at a point $Q \in Y$. While $Y$ is not NPC, $T_QY$ is always NPC by construction.

To prove item (2) in Theorem \ref{LipThm}, we first note that all of the results of \cite[section 6]{daskal-meseCAG} can be immediately applied. In particular, set $f_k:= f_{\lambda_k}$ and define $h_k :({\bf B}(1),\delta) \to (\mathcal CY, D_{\lambda_k})$ to be the minimizer into $\mathcal CY$ with $h_k|_{\partial {\bf B}(1)}= \overline f_k|_{\partial {\bf B}(1)}$. Then $h_k$ converges uniformly locally in the pull-back sense to $\overline f_*$ and in fact the pseudodistance functions $d_*^f, d_*^h$, determined by the maps $f_k$ and $h_k$, are equal.

The results of \cite[section 8.3]{daskal-meseCAG} rely on the fact that the sequence $f_k$ satisfies the monotonicity formula and that the tangent map $f_*$ is homogeneous. We have already established monotonicity for the sequence $f_k$. To prove homogeneity for the tangent map $\overline f_*$, we proceed exactly as in \cite{daskal-meseCAG}. 

\begin{lemma}\label{tangentenergy}
For ${\bf B}$ a dimension-$n$, codimension-$\nu$ local model, $\nu \geq 2$, and $g$ a normalized Lipschitz metric defined on ${\bf B}(r)$, $(Y,d)$ a CAT(1) space, let $f:({\bf B}(r),g) \to \B_\tau(P) \subset Y$ be an energy minimizing map where $0<\tau< \frac \pi 4$. Let $f_k, h_k,\overline f_*$ be as above. Then the directional energies of the sequences $h_k, f_k$ converge to the directional energies of $\overline f_*$. 
\end{lemma}
\begin{proof}
Since the $h_k$ are H\"older continuous and satisfy a monotonicity formula, the proof of \cite[Lemma 8.8]{daskal-meseCAG} can be followed verbatim. To prove the directional energies converge relies only on estimates relating the energies of $h_k, f_k$, and lower semi-continuity of the energy. Since the energy comparisons follow from the comments above, the result is immediate.
\end{proof}
Now, following the proof of \cite[Lemma 6.3]{daskal-meseCAG}, since the directional energies of $h_k$ converge to those of $\overline f_*$, $\overline f_*$ is homogeneous of order $\alpha$, i.e.
\[
d_*(\overline f_*(x), \overline f_*(0)) = |x|^\alpha d_*\left(\overline f_*\left(\frac x{|x|}\right), \overline f_*(0)\right).
\]
 Note that in the proof in \cite{daskal-meseCAG}, the right hand side of the equation for $(E^{f_*}(\sigma))'$ should include the term $\frac{n-2}\sigma E^{f_*}(\sigma)$ and the first term in the parenthesis in (6.10) should be the product of boundary integrals.
 
With the homogeneity in hand, we can now follow the proofs of \cite{daskal-meseCAG} to conclude the necessary results. For ${\bf B}$ a dimension-$n$, codimension-$\nu$ local model, recall that $D$ is isomorphic to  $\mathbb R^{n-\nu}$. For each $x \in D$, let $N(x)$ be the $\nu$-plane orthogonal to $D$ at $x$. Then, for $|x|<1$, $\partial {\bf B}(1) \cap N(x)$ is a spherical $(\nu-1)$-complex. The key proposition, which follows exactly the proof of \cite[Theorem 8.4]{daskal-meseCAG}, is as follows. 
\begin{proposition}
Let ${\bf B}$ be a dimension-$n$, codimension-$\nu$ local model, $\nu \geq 2$, and $g$ a normalized Lipschitz metric defined on ${\bf B}(r)$ and $(Y,d)$ a CAT(1) space, let $f:({\bf B}(r),g) \to \B_\tau(P) \subset Y$ be an energy minimizing map where $0<\tau< \frac \pi 4$. If $\lambda_1(\partial {\bf B}(1) \cap N(0),T_QY) \geq \beta (>\beta)$ for all $Q \in Y$ and $\alpha <1$, then the order $\alpha$ of $f$ at $0$ satisfies $\alpha(\alpha+\nu-2) \geq \beta (>\beta)$. 
\end{proposition} 

 For a local model ${\bf B}$ and any $x \in {\bf B}(r)$, recall that $f \circ L_x:{\bf B}_x'(r(x)) \to \B_\tau P$ where $L_x$ is given by \cite[Proposition 2.1]{daskal-meseCAG} and ${\bf B}_x'(r(x))$ is a local model centered at $x$. Define
\[
\lambda_1^N:= \inf_{x \in {\bf B}(r) \cap N(0), Q \in Y} \lambda_1(\partial {\bf B}_x'(1), T_QY).
\]As an immediately corollary of the previous proposition, we observe that:
\begin{corollary}
For ${\bf B}$ a dimension-$n$, codimension-$\nu$ local model, $\nu \geq 2$, $g$ a Lipschitz Riemannian metric defined on ${\bf B}(r)$, and $(Y,d)$ a CAT(1) space, let $f:({\bf B}(r),g) \to \B_\tau(P) \subset Y$ be an energy minimizing map where $0<\tau< \frac \pi 4$. If $\lambda_1^N \geq \nu -1$ then $f$ is Lipschitz continuous in ${\bf B}(\varrho r)$ for $\varrho \in (0,1)$.
\end{corollary}

Theorem \ref{LipThm}, item (2), immediately follow from the above results following the observation that $\nu = n-k$.


\end{document}